\documentclass[11pt,a4paper,leqno]{amsart}
       
\usepackage{amsmath,amsfonts,amssymb}
\usepackage{enumerate} 
\usepackage{amsthm}
\usepackage{hyperref}

\def\R{{\mathbb R}}
\def\N{{\mathbb N}}

\def\virgp{\raise 2pt\hbox{,}}
\def\bv{{\bf v}}
\def\bu{{\bf u}}
\def\bw{{\bf w}}
\def\({\left(}
\def\){\right)}
\def\<{\left\langle}
\def\>{\right\rangle}

\def\Eq#1#2{\mathop{\sim}\limits_{#1\rightarrow#2}}

\def\d{{\partial}}
\def\a{{\tt a}}

\def\e{\varepsilon}

\def\om{\omega}

\def\F{\mathcal F}

\def\O{\mathcal O}
\def\eik{\phi_{\rm eik}}

\theoremstyle{plain}
\newtheorem{theo}{Theorem}[section]
\newtheorem{lem}[theo]{Lemma}
\newtheorem{cor}[theo]{Corollary}
\newtheorem{prop}[theo]{Proposition}
\newtheorem{hyp}{Assumption}

\theoremstyle{definition}

\theoremstyle{remark}
\newtheorem{rema}[theo]{Remark}
\newtheorem*{rema*}{Remark}

\newtheorem{ex}{Example}

\numberwithin{equation}{section}

\begin{document}
\author[R. Carles]{R{\'e}mi Carles}
\address{MAB, Universit\'e de Bordeaux 1\\ 351 cours de la
  Lib{\'e}ration\\ 33405 Talence cedex\\ France}
\email{Remi.Carles@math.cnrs.fr}
\thanks{Supports by European network HYKE,
  funded  by the EC as contract HPRN-CT-2002-00282, by Centro de
  Matem\'atica e Aplica\c c\~oes Fundamentais (Lisbon), funded by
  FCT as contract POCTI-ISFL-1-209, and by the ANR project SCASEN,
  are acknowledged.} 
\title[WKB for NLS]{WKB analysis for nonlinear Schr\"odinger equations
with potential} 
\begin{abstract}
We justify the WKB analysis for the semiclassical nonlinear
Schr\"o\-din\-ger equation with a subquadratic potential. This concerns
subcritical, critical, and supercritical cases as far as the
geometrical optics method is concerned. In the supercritical case,
this extends a previous result by E.~Grenier; we also have to
restrict to nonlinearities which are defocusing and cubic at the
origin, but besides subquadratic potentials, we consider
initial phases which may be unbounded. For this, we construct solutions
for some compressible Euler equations with unbounded 
source term and unbounded initial velocity.   
\end{abstract}
\maketitle

\numberwithin{equation}{section}

\section{Introduction}
\label{sec:intro}
Consider the initial value problem, for $x\in \R^n$ and $\kappa \ge 0$:
\begin{align}
  i\e \d_t u^\e +\frac{\e^2}{2}\Delta u^\e &= V(t,x)u^\e + \e^\kappa
  f\(|u^\e|^2\) 
  u^\e\label{eq:nlssemi}\\
u^\e_{\mid t=0} &= a_0^\e(x)e^{i\phi_0(x)/\e} .\label{eq:CI} 
\end{align}
The aim of WKB methods is to describe $u^\e$ in the limit $\e \to 0$,
when $\phi_0$ does not depend on $\e$, and 
$a_0^\e$ has an asymptotic expansion of the form:
\begin{equation}\label{eq:CIbkw}
  a_0^\e (x)\sim a_0(x) +\e a_1(x)+\e^2 a_2(x) +\ldots
\end{equation}
The parameter $\kappa \ge 0$ describes the strength of a coupling
constant, which makes nonlinear effects more or less important in the
limit $\e \to 0$; the larger the $\kappa$, the weaker the nonlinear
interactions. In this paper, we describe the asymptotic behavior of
$u^\e$ at leading order, when the potential $V$ and the initial phase
$\phi_0$ are smooth, and subquadratic in the space variable. 

\smallbreak

Such equations as \eqref{eq:nlssemi} appear in physics: see
e.g. \cite{Sulem} for a general overview.  For instance, they are used to
model Bose-Einstein condensation when $V$ is an
harmonic potential (isotropic or anisotropic) and the nonlinearity is
cubic or quintic (see e.g. \cite{DGPS,KNSQ,PiSt}). In most of this
paper, the initial data
we consider are in Sobolev spaces $H^s$.  
We refer to \cite{BJM} for numerics on the semi-classical limit of
\eqref{eq:nlssemi}. 

\smallbreak

We shall not recall the results concerning the  Cauchy
problem for \eqref{eq:nlssemi}-\eqref{eq:CI}, and refer to
\cite{CazCourant} 
for an overview on the semilinear Schr\"odinger equation. 

\smallbreak

In the one-dimensional case $n=1$, the cubic nonlinear
Schr\"odinger equation is completely integrable, in the absence of an
external potential,  or when $V$ is a quadratic polynomial
(\cite[p.~375]{AblowitzClarkson}). Several tough papers
analyze the semi-classical limit in the case $V\equiv 0$, for $\kappa =0$: see
e.g. \cite{JLM,KMM,TianYe,TVZ}. We shall not use this approach, but
rather work in the spirit of \cite{Grenier98}.  

\smallbreak

An interesting feature of \eqref{eq:nlssemi} is that one does not
expect the creation of harmonics, provided that only one phase is
present initially, like in \eqref{eq:CI}. The WKB methods consist in
seeking an approximate solution to \eqref{eq:nlssemi} of the form:
\begin{equation}\label{eq:defBKW}
  u^\e(t,x)\sim \(\a_0(t,x) + \e \a_1(t,x) + \e^2
  \a_2(t,x)+\ldots \) e^{i\Phi(t,x)/\e}\, .
\end{equation}
One must not expect this approach to be valid when caustics are
formed: near a caustic, all the terms
$\Phi$, $\a_0$, $\a_1$, \ldots become singular. Past the
caustic, several phases are necessary in general to describe the
asymptotic behavior of the solution (see e.g. \cite{Du} for a general
theory in the linear case). In this paper,
we restrict our attention to times preceding this break-up. 
\smallbreak

For such an expansion to be available with profiles $\a_j$
independent of $\e$, it is reasonable to assume that $\kappa$ is an
integer. However, we do not assume that $\kappa$ is an integer:
we study the asymptotic behavior of $u^\e$ \emph{at leading order}
(strong limits in $L^2\cap 
L^\infty$ for instance), including cases where other powers of $\e$
would come into play.  
\smallbreak

We distinguish two families of assumptions: ``geometrical''
assumptions, on the potential $V$ and the initial phase $\phi_0$, and
``analytical'' assumptions, on $f$ and the initial amplitude
$a_0^\e$. In all the cases, we shall not try to seek the optimal
regularity; we focus our interest on the limit $\e\to 0$.

\begin{hyp}[Geometrical]\label{hyp:geom}
  We assume that the potential and the initial phase are smooth and
  subquadratic: 
  \begin{itemize}
  \item $V\in C^\infty(\R_t\times \R^n_x)$, and $\d_x^\alpha V\in
  L^\infty_{\rm loc}(\R_t ;L^\infty(\R^n_x))$ as soon as $|\alpha|\ge 2$.
  \item $\phi_0\in C^\infty( \R^n)$, and $\d_x^\alpha \phi_0\in
  L^\infty(\R^n)$ as soon as $|\alpha|\ge 2$.
  \end{itemize}
\end{hyp}
The assumption of $V$ being subquadratic is classical in other
contexts; for instance, locally in time, the dispersion for
$e^{-i\frac{t}{\e}(-\e^2\Delta +V)}$ is the same as without potential
(see \cite{Fujiwara79,Fujiwara}), 
\begin{equation*}
  \left\| e^{-i\frac{t}{\e}(-\e^2\Delta +V)}\right\|_{L^1\to L^\infty}
  \le \frac{C}{|\e t|^{n/2}},\quad \forall |t|\le \delta,
\end{equation*}
hence yielding the same local Strichartz
estimates as in the free case. Global in time Strichartz
  estimates must not be expected in general, as shown by the example
  of the harmonic oscillator, which has eigenvalues.  For
  positive superquadratic potentials, the 
smoothing effects and Strichartz estimates are different (see
\cite{YajZha01,YajZha04}). This is related to the 
properties of the Hamilton flow, which also imply: 
\begin{lem}\label{lem:hj}
  Under Assumption~\ref{hyp:geom}, there exist $T>0$ and a unique
  solution $\eik\in C^\infty([0,T]\times\R^n)$ to:
  \begin{equation}
    \label{eq:eik}
    \d_t \eik +\frac{1}{2}|\nabla_x \eik|^2 +V(t,x)=0\quad ;\quad
    \phi_{{\rm eik} \mid t=0}=\phi_0\, .
  \end{equation}
This solution is subquadratic: $\d_x^\alpha \eik \in
  L^\infty([0,T]\times\R^n)$ as soon as $|\alpha|\ge 2$. 
\end{lem}
This result is proved in Section~\ref{sec:HJ}, where other remarks on
Assumption~\ref{hyp:geom} are made. 

\begin{hyp}[Analytical]\label{hyp:ana}
  We assume that the nonlinearity is smooth,  and
  that the initial amplitude converges in
  Sobolev spaces: 
  \begin{itemize}
  \item $f\in C^\infty(\R; \R)$.
  \item There exists  $ a_0\in H^\infty := \cap_{s\ge
  0}H^s(\R^n)$, such that 
  $ a_0^\e$ converges to  $a_0$ in $H^s$ for any $s$, as
  $\e\to 0$.  
  \end{itemize}
\end{hyp}
\begin{rema*}
  Some of the results we shall prove remain valid when $f$ is
  complex-valued. In that case, the conservation of mass associated to
  the Schr\"odinger equation, 
  $\|u^\e(t)\|_{L^2}= \|a_0^\e\|_{L^2}$, no longer holds. On the other
  hand, when $0\le \kappa <1$, this assumption is necessary in
  our approach, and we even assume $f'>0$. 
\end{rema*}
\subsection{Subcritical and critical cases: $\kappa \ge 1$}
\label{sec:sub0}

When the initial data is of the form \eqref{eq:CIbkw}, the usual
approach consists in plugging a formal expansion of the form
\eqref{eq:defBKW} into \eqref{eq:nlssemi}. Ordering the terms in
powers of $\e$, and canceling the cascade of equations thus obtained
yields $\Phi$, $\a_0$, $\a_1$, \ldots 
\smallbreak

Assume in this section
that $\kappa \ge 1$, and apply the above procedure. To cancel the term
of order $\O(\e^0)$, we find that $\Phi$ must solve
\eqref{eq:eik}: $\Phi=\eik$. Canceling the term of order $\O(\e^1)$,
we get: 
\begin{equation*}
  \d_t \a_0 +\nabla \eik\cdot \nabla \a_0 +
  \frac{1}{2}\a_0\Delta \eik =
\left\{
  \begin{aligned}
    &0 & \text{ if }\kappa >1,\\
    &-if\(|\a_0|^2\)\a_0& \text{ if }\kappa =1.
  \end{aligned}
\right.
\end{equation*}
We see that  the value $\kappa =1$ is critical as far as
nonlinear effects are concerned: if $\kappa>1$, no nonlinear effect is
expected at leading order, since formally, $u^\e \sim \a_0 e^{i\eik
  /\e}$, and $\eik$ and $\a_0$ do not depend on the nonlinearity
$f$. If $\kappa =1$, then $\a_0$ solves a nonlinear equation
involving $f$. 
\smallbreak

We will see in Section~\ref{sec:sub} that $\a_0$ solves a transport
equation that turns out to be a ordinary differential equation along
the rays of geometrical optics, as is usual in the hyperbolic case
(see e.g. \cite{RauchUtah}). More typical of Schr\"odinger equation is
the fact that this ordinary differential equation can be solved
explicitly: the nonlinear effect is measured by a nonlinear phase
shift (see the example of \cite{LannesRauch} for a similar result 
in the hyperbolic setting). We prove the following result in
Section~\ref{sec:sub}: 
\begin{prop}\label{prop:sub}
  Let Assumptions~\ref{hyp:geom} and \ref{hyp:ana} be satisfied. Let
  $\kappa \ge 1$. Then for any $\e \in ]0,1]$,
  \eqref{eq:nlssemi}-\eqref{eq:CI} has a unique solution $u^\e \in
  C^\infty([0,T]\times \R^n)\cap C([0,T];H^s)$ for any $s>n/2$ ($T$ is
  given by  
  Lemma~\ref{lem:hj}). Moreover,
  there exist $ a, G\in C^\infty([0,T]\times \R^n)$,
  independent of $\e \in ]0,1]$, where
  $ a \in C([0,T]; L^2\cap L^\infty)$, and $G$ is
  real-valued with 
  $G \in C([0,T]; 
  L^\infty)$, such that:
  \begin{equation*}
    \left\| u^\e -  a e^{i\e^{\kappa -1}G}e^{i\eik
        /\e}\right\|_{L^\infty([0,T]; L^2\cap L^\infty) } \to 0\quad
    \text{as }\e \to 0.
  \end{equation*}
The profile $ a$ solves the initial value problem:
\begin{equation}\label{eq:alibre}
  \d_t  a +\nabla \eik\cdot \nabla a +
  \frac{1}{2} a\Delta \eik =0\quad ;\quad a_{\mid
    t=0}=a_0, 
\end{equation}
and $G$  depends nonlinearly on $ a$ through $f$.
In particular, if $\kappa>1$, then 
\begin{equation*}
    \left\| u^\e - a e^{i\eik
        /\e}\right\|_{L^\infty([0,T]; L^2\cap L^\infty) } \to 0\quad
    \text{as }\e \to 0,
  \end{equation*}
and no nonlinear effect is present in the leading order behavior of
$u^\e$. If $\kappa =1$, nonlinear effects are present at leading
order, measured by $G$.
\end{prop}
The dependence of $G$ upon $ a$ and $f$ is made more explicit in
Section~\ref{sec:sub}, in terms of the Hamilton flow determining
$\eik$ (see \eqref{eq:exprG}). Note that in the above result, we do
not assume that $\kappa$ is an integer. 

\subsection{Supercritical case: $\kappa=0$}
\label{sec:surcrit0}

It follows from the above analysis that the case $0\le \kappa <1$ is
supercritical. We restrict our attention to the case $\kappa =0$. We
present an analysis of the range $0<\kappa<1$ in
Section~\ref{sec:kappaint}. Plugging an
asymptotic expansion of the form \eqref{eq:defBKW} into
\eqref{eq:nlssemi} yields a shifted cascade of equations:
\begin{align*}
 \O\(\e^0\):&\quad \d_t \Phi +\frac{1}{2}|\nabla\Phi|^2 + V +
 f\(|\a_0|^2\)=0,\\ 
\O\(\e^1\):&\quad \d_t \a_0 +\nabla\Phi \cdot \nabla \a_0
+\frac{1}{2}\a_0\Delta \Phi = 
2if'\(|\a_0|^2\) {\rm Re}\(\a_0\overline{\a_1}\).
\end{align*}
Two comments are in order. First, we see that there is a strong
coupling between the phase and the main amplitude: $\a_0$ is present
in the equation for $\Phi$. Second, the above system is not closed:
$\Phi$ is determined in function of $\a_0$, and $\a_0$ is determined in
function of $\a_1$. Even if we pursued the cascade of equations, this
phenomenon would remain: no matter how many terms are computed, the
system is never closed (see \cite{PGX93}). This is a typical feature
of supercritical 
cases in nonlinear geometrical optics (see
\cite{CheverryBullSMF,CG05}).  
\smallbreak

In the case when $V\equiv 0$ and $\phi_0\in H^s$, this problem was
resolved by E.~Grenier \cite{Grenier98}, by modifying the usual WKB
methods; this approach is recalled in Section~\ref{sec:grenier}. Note
that even though $\a_1$ is not determined by the above system, the
pair $(\rho,v):= (|\a_0|^2,\nabla \Phi)$ solves a compressible Euler
equation:
\begin{equation}
  \label{eq:eulerV}
  \begin{aligned}
    \d_t v +v\cdot \nabla v + \nabla V + \nabla f(\rho)&=0\ ;\quad
    v\big|_{ t=0}=\nabla \phi_0\\
    \d_t \rho+ \nabla\cdot (\rho v) &=0\ ;\quad \rho\big|_{t=0}=|a_0|^2.
  \end{aligned}
\end{equation}
Using techniques introduced in the study of quasilinear hyperbolic
equations, E.~Grenier justified a WKB expansion for nonlinearities
which are defocusing, and cubic at the origin ($f'>0$). We shall not
change this assumption, but show how to treat the case of a
subquadratic potential with a subquadratic initial phase. Note that
even the construction of solution to \eqref{eq:eulerV} under
Assumption~\ref{hyp:geom} is not standard: the source term $ \nabla V$
may be unbounded, as well as the initial velocity $\nabla \phi_0$. 
\begin{hyp}\label{hyp:surcrit}
  In addition to Assumption~\ref{hyp:ana}, we assume:
  \begin{itemize}
  \item $f'>0$.
  \item There exists $a_0, a_1\in H^\infty$, with $xa_0, xa_1\in H^\infty$,
  such that:
  \begin{equation*}
    \left\|a_0^\e -a_0 -\e a_1\right\|_{H^s} +
    \left\|xa_0^\e -xa_0 -\e xa_1\right\|_{H^s}
    =o(\e),\quad \forall s\ge 0. 
  \end{equation*}
  \end{itemize}
\end{hyp}
We can then describe the asymptotic behavior of the solution to
\eqref{eq:nlssemi}-\eqref{eq:CI} for small times:
\begin{theo}\label{theo:BKWV}
  Let Assumptions~\ref{hyp:geom}, \ref{hyp:ana} and \ref{hyp:surcrit}
  be satisfied. Let $\kappa =0$. There exists $T_*>0$ independent
  of $\e \in ]0,1]$ and a unique solution $u^\e \in
  C^\infty([0,T_*]\times \R^n)\cap C([0,T_*];H^s)$ for any $s>n/2$ to
  \eqref{eq:nlssemi}-\eqref{eq:CI}. Moreover, there exist $a,\phi \in
  C([0,T_*];H^s)$ for every $s\ge 0$, such that:
  \begin{equation}\label{eq:BKWVtpetit}
    \limsup_{\e \to 0}\left\| u^\e - a
    e^{i(\phi+\eik)/\e}\right\|_{L^2\cap L^\infty}=\O(t)\quad
    \text{as }t\to 0. 
  \end{equation}
Here, $a$ and $\phi$ are nonlinear functions of $\eik$ and
$a_0$, given by \eqref{eq:11h09}. Finally, there exists $\phi^{(1)}\in
  C([0,T_*];H^s)$ for every $s\ge 0$, real-valued, such that: 
  \begin{equation}\label{eq:BKWVtgrand}
    \limsup_{\e \to 0}\sup_{0\le t\le T_*}\left\| u^\e - ae^{i\phi^{(1)}}
    e^{i(\phi+\eik)/\e}\right\|_{L^2\cap L^\infty}=0. 
  \end{equation}
The phase shift $\phi^{(1)}$ is a nonlinear function of $\eik, a_0$
and $a_1$.  
\end{theo}
This result can be understood as follows. At leading order, the
amplitude of $u^\e$ is given by $a e^{i\phi^{(1)}}$, which can be
approximated for small times by $a$, because
$\phi^{(1)}\big|_{t=0}=0$. The rapid oscillations are described by the
phase $\phi +\eik$. The 
function $\phi$ is constructed as a perturbation of $\eik$, but must not
be considered as negligible: its $H^s$-norms are not small in
general, see \eqref{eq:11h09} (at time $t=0$ for instance). As a
consequence of our analysis, the pair 
\begin{equation*}
  (\rho,v) = \( |a|^2,\nabla (\phi +\eik)\) = \(
  \left|a e^{i\phi^{(1)}}\right|^2,\nabla (\phi +\eik)\)  
\end{equation*}
solves the system \eqref{eq:eulerV}. 
\begin{rema}
With this result, we could deduce instability phenomena for 
\eqref{eq:nlssemi}-\eqref{eq:CI} in the same fashion as in
\cite{CaInstab}. Note that because of Assumption~\ref{hyp:geom}, it 
seems that the
approaches of \cite{BGTENS,BZ,CCT2} cannot be adapted to the
present case: the Laplacian can never be neglected, and apparently, 
WKB approach is really needed.  
\end{rema}

The rest of this paper is organized as follows. In
Section~\ref{sec:HJ}, we prove Lemma~\ref{lem:hj}. In
Section~\ref{sec:sub}, we prove Proposition~\ref{prop:sub}, and
explain  how $G$ is obtained. We recall the approach of
\cite{Grenier98} in Section~\ref{sec:grenier}, and show how to adapt
it to prove Theorem~\ref{theo:BKWV} 
in Section~\ref{sec:surcrit}. 
We present an
analysis for the case $0<\kappa<1$ in Section~\ref{sec:kappaint}. 


\section{Global in space Hamilton-Jacobi theory}
\label{sec:HJ}
In this section, we prove Lemma~\ref{lem:hj}.
Consider the classical Hamiltonian:
\begin{equation*}
  H(t,x,\tau,\xi)= \tau +\frac{1}{2}|\xi|^2 +V(t,x),\quad
  (t,x,\tau,\xi)\in \R_+\times\R^n \times\R\times\R^n.
\end{equation*}
It is smooth by Assumption~\ref{hyp:geom}. Therefore, it is classical
(see e.g. \cite{DG}) that in the neighborhood of each
point $x\in \R^n$, one can construct a smooth solution to the eikonal
equation \eqref{eq:eik}, on some time interval $[-t(x),t(x)]$, for
some $t(x)>0$ depending on $x$. The fact 
that in Lemma~\ref{lem:hj}, we can find some $T>0$ uniform in $x\in
\R^n$ is due to the fact that the potential and the initial phase are
subquadratic. 
\smallbreak

Recall that if there
exist some constants $a,b>0$ such that a potential $\tt V$ satisfies
${\tt V}(x)\geq -a |x|^2-b$, then $-\Delta +{\tt V}$ is
essentially self-adjoint on $C_0^\infty({\mathbb R}^n)$ (see 
\cite[p.~199]{ReedSimon2}). If  $-{\tt V}$ has
superquadratic growth,
then it is not possible to define $e^{-it(-\Delta +{\tt V})}$
(see \cite[Chap.~13, Sect.~6, Cor.~22]{Dunford} for the case ${\tt
  V}(x)=-x^4$ in space dimension one). This is due to the fact that
classical trajectories can reach an infinite speed. We will see below
that if the initial phase $\phi_0$ is superquadratic, then focusing at
the origin may occur ``instantly'' (see Example~\ref{ex:phase}). 
\smallbreak

To construct the solution of \eqref{eq:eik}, introduce the flow
associated to $H$: let $x(t,y)$ and $\xi(t,y)$ solve
\begin{equation}
  \label{eq:hamilton}
\left\{
  \begin{aligned}
   &\d_t x(t,y) = \xi \(t,y\)\quad ;\quad x(0,y)=y,\\ 
   &\d_t \xi(t,y) = -\nabla_x V\(t,x(t,y)\)\quad ;\quad
   \xi(0,y)=\nabla \phi_0(y).
  \end{aligned}
\right.
\end{equation}
Recall the result (valid under weaker conditions than
Assumption~\ref{hyp:geom}): 
\begin{theo}[\cite{DG}, Th.~A.3.2]\label{theo:hj}
  Suppose that Assumption~\ref{hyp:geom} is satisfied. Let $t\in
  [0,T]$ and $\theta_0$ an open set of $\R^n$. Denote
  \begin{equation*}
    \theta_t := \{ x(t,y)\in \R^n, y\in\theta_0\}\quad ;\quad \theta
    := \{ (t,x)\in [0,T]\times\R^n, x\in \theta_t\}.
  \end{equation*}
Suppose that for $t\in [0,T]$, the mapping 
\begin{equation*}
  \theta_0\ni y\mapsto x(t,y)\in \theta_t
\end{equation*}
is bijective, and denote by $y(t,x)$ its inverse. Assume also that 
\begin{equation*}
  \nabla_x y \in L^\infty_{\rm loc}(\theta).
\end{equation*}
Then there exists a unique function $\theta\ni (t,x)\mapsto \eik(t,x)\in
\R$ that solves \eqref{eq:eik}, and satisfies $\nabla_x^2\eik \in
L^\infty_{\rm loc}(\theta)$. Moreover,
\begin{equation}\label{eq:gradeik}
  \nabla_x \eik(t,x) =\xi(t,y(t,x)). 
\end{equation}
\end{theo}

\begin{prop}[\cite{SchwartzBook}, Th.~1.22 and \cite{DG},
  Prop.~A.7.1]\label{prop:invglobale} 
  Suppose that the function $\R^n\ni y\mapsto x(y)\in \R^n$ satisfies:
  \begin{equation*}
    |\operatorname{det}\nabla_y x|\ge C_0>0\quad \text{and}\quad \left|
     \d_y^\alpha x\right|\le C, \ |\alpha|=1,2.
  \end{equation*}
Then $x$ is bijective. 
\end{prop}
\begin{proof}[Proof of Lemma~\ref{lem:hj}]
Lemma~\ref{lem:hj} follows from the above two
results. From Assumption~\ref{hyp:geom}, we know that we can solve
\eqref{eq:hamilton} locally in time in the neighborhood of any
$y\in\R^n$. Differentiate \eqref{eq:hamilton} with respect to $y$:
\begin{equation}
  \label{eq:dyhamilton}
\left\{
  \begin{aligned}
   &\d_t \d_y x(t,y) = \d_y\xi \(t,y\)\quad ;\quad \d_y
   x(0,y)={\rm Id},\\ 
   &\d_t \d_y \xi(t,y) = -\nabla_x^2 V\(t,x(t,y)\)\d_y x(t,y) \quad ;\quad
   \d_y\xi(0,y)=\nabla^2 \phi_0(y).
  \end{aligned}
\right.
\end{equation}
Integrating \eqref{eq:dyhamilton} in time, we infer from
Assumption~\ref{hyp:geom} that for any $T>0$, there exists $C_T$ such
that for $(t,y)\in [0,T]\times \R^n$: 
\begin{equation*}
  \left| \d_y x(t,y) \right| +\left|  \d_y \xi(t,y)\right|
  \le C_T + C_T\int_0^t \( \left| \d_y x(s,y) \right| +\left|  \d_y
  \xi(s,y)\right| \) ds.
\end{equation*}
Gronwall lemma yields:
\begin{equation}\label{eq:10:04}
  \left\| \d_y x(t) \right\|_{L^\infty_y} +\left\|  \d_y
  \xi(t)\right\|_{L^\infty_y} 
  \le C'(T).
\end{equation}
Similarly, 
\begin{equation}\label{eq:10:05}
  \left\| \d_y^\alpha x(t) \right\|_{L^\infty_y} +\left\|
  \d_y^\alpha  \xi(t)\right\|_{L^\infty_y} 
  \le C(\alpha,T),\quad \forall \alpha \in \N^n,\  |\alpha|\ge 1.
\end{equation}
Integrating the first line of \eqref{eq:dyhamilton} in time, we
have:
\begin{equation*}
  \operatorname{det}\nabla_y x(t,y) = \operatorname{det}\({\rm
  Id}+\int_0^t \nabla_y\xi \(s,y\)ds\). 
\end{equation*}
We infer from \eqref{eq:10:04} that for $t\in [0,T]$, provided that
$T>0$ is sufficiently small, we can find $C_0>0$ such that:
\begin{equation}\label{eq:det}
  \left| \operatorname{det}\nabla_y x(t,y)\right|\ge C_0,\quad \forall
  (t,y)\in [0,T]\times \R^n.
\end{equation}
Applying Proposition~\ref{prop:invglobale}, we deduce that we can
invert $y\mapsto x(t,y)$ for $t\in [0,T]$. 
\smallbreak

To apply Theorem~\ref{theo:hj} with $\theta_0= \theta= \theta_t=\R^n$,
we must check that $\nabla_x y\in L^\infty_{\rm loc}
(\R^n)$. Differentiate the relation
\begin{equation*}
  x\(t,y(t,x)\) = x
\end{equation*}
with respect to $x$:
\begin{equation*}
  \nabla_x y(t,x)\nabla_y x\(t,y(t,x)\) = {\rm Id}. 
\end{equation*}
Therefore, $\nabla_x y(t,x) = \nabla_y x\(t,y(t,x)\)^{-1}$ as
matrices, and
\begin{equation}\label{eq:adj}
  \nabla_x y(t,x) =\frac{1}{\operatorname{det}\nabla_y x(t,y)}{\rm
  adj}\( \nabla_y x\(t,y(t,x)\)\),
\end{equation}
where ${\rm adj}\( \nabla_y x\)$ denotes the adjugate of $\nabla_y
x$. We infer from \eqref{eq:10:04} and \eqref{eq:det} that $\nabla_x y
\in L^\infty(\R^n)$ for $t\in [0,T]$. Therefore, Theorem~\ref{theo:hj}
yields a smooth solution $\eik$ to \eqref{eq:eik}, local in time and
global in space: $\eik \in C^\infty([0,T]\times \R^n)$. 
\smallbreak

The fact that $\eik$ is subquadratic as stated in Lemma~\ref{lem:hj}
then stems from \eqref{eq:gradeik}, \eqref{eq:10:05}, \eqref{eq:det}
and \eqref{eq:adj}. 
\end{proof}

We now give some two examples showing that
Assumption~\ref{hyp:geom} is essentially sharp to solve \eqref{eq:eik}
globally in space, at least when no assumption is made on the sign of
$V$ nor on the geometry of $\nabla \phi_0$. We already recalled that
if $-V$ has a 
superquadratic growth, then $-\Delta +V$ is not essentially
self-adjoint on $C_0^\infty(\R^n)$, so we shall rather study the
dependence of $\eik$ on the initial phase $\phi_0$.
\begin{ex}\label{ex:phase}
  Assume that $V\equiv 0$ and 
  \begin{equation*}
    \phi_0(x) = -\frac{1}{(2+2\delta)T}\(|x|^2+1 \)^{1+\delta}, \quad
    T>0,\ \delta
    \not = -1. 
  \end{equation*}
Then Assumption~\ref{hyp:geom} is satisfied if and only if $\delta\le
0$. When $\delta =0$, then \eqref{eq:eik} is solved explicitly:
\begin{equation*}
  \eik(t,x) = \frac{|x|^2}{2(t-T)}-\frac{1}{2T}. 
\end{equation*}
This shows that we can solve \eqref{eq:eik} globally in space, but only
locally in time: as $t\to T$, a caustic reduced to a single point (the
origin) is
formed. Note that with $T<0$, \eqref{eq:eik} can be solved globally in
time \emph{for positive times}.\\
When $\delta>0$, then integrating \eqref{eq:hamilton} yields:
\begin{align*}
 x(t,y) &=y +\int_0^t \xi(s,y)ds = y +\int_0^t  \xi(0,y)ds= y -
\frac{t}{T}\(|y|^2+1 \)^{\delta} y\\
&= y\(1 - \frac{t}{T}\(|y|^2+1
\)^{\delta}\).  
\end{align*}
For $R>0$, we see that the rays starting from the
ball $\{|y|=R\}$ meet at the origin at time
\begin{align*}
  T_c(R) = \frac{T}{(R^2+1)^\delta}.
\end{align*}
Since $R$ is arbitrary, this shows that several rays can meet
arbitrarily fast, thus showing that Theorem~\ref{theo:hj} cannot be
applied uniformly in space. 
\end{ex}
\begin{ex}
  When $V(t,x) = \frac{1}{2}\sum_{j=1}^n \om_j^2 x_j^2$ is an harmonic
  potential ($\om_j \not= 0$), and $\phi_0\equiv 0$, we have:
  \begin{equation*}
    \phi(t,x) = -\sum_{j=1}^n \frac{\om_j}{2}x_j^2\tan (\om_j t). 
  \end{equation*}
This also shows that we can solve \eqref{eq:eik} globally in space,
but locally in time only. Note that if we replace formally $\om_j$ by
$i\om_j$, then $V$ is turned into $-V$, and the trigonometric
functions become hyperbolic functions: we can then solve \eqref{eq:eik}
globally in space \emph{and} time. 
\end{ex}

Instead of invoking Theorem~\ref{theo:hj} and
Proposition~\ref{prop:invglobale}, one might try to differentiate
\eqref{eq:eik} in order to prove that $\eik$ is subquadratic, in the
same fashion as in \cite{BahouriCheminAJM,BahouriCheminIMRN}. For
$1\le j,k\le n$, differentiate \eqref{eq:eik} with respect to $x_j$
and $x_k$:
\begin{align*}
    &\d_t \d^2_{jk}\eik +\nabla_x \eik\cdot \nabla_x\(\d^2_{jk}\eik \)+
    \sum_{l=1}^n \d^2_{jl}\eik\d^2_{lk}\eik
    +\d^2_{jk}V(t,x)=0\ ;\\
    &\d^2_{jk}\phi_{{\rm eik} \mid t=0}=\d^2_{jk}\phi_0\, .
  \end{align*}
We see that we obtain a system of the form
\begin{equation*}
  D_t y = Q(y) + R\quad ;\quad y_{\mid t=0} = y(0),
\end{equation*}
where $y$ stands for the family $(\d^2_{jk}\eik)_{1\le j,k\le n}$, $Q$
is quadratic, and
$R$ and $y(0)$ are bounded. The operator $D_t$ is a well-defined
transport operator provided that the characteristics given by:
\begin{equation*}
  \d_t x(t,y) = \nabla_x \eik \(t,x(t,y)\)\quad ; \quad x(0,y)=y,
\end{equation*}
define a global diffeomorphism. Proving this amounts to using
Proposition~\ref{prop:invglobale}, for the rather general initial data
we consider. So it seems that this approach does not allow to shorten
the proof of Lemma~\ref{lem:hj}.


\section{Subcritical and critical cases}
\label{sec:sub}

To establish Proposition~\ref{prop:sub}, define
\begin{equation*}
  a^\e(t,x) := u^\e (t,x) e^{-i\eik(t,x)/\e}. 
\end{equation*}
Then $u^\e$ solves \eqref{eq:nlssemi}-\eqref{eq:CI} if and only if
$a^\e$ solves:
\begin{equation}
  \label{eq:asub}
  \begin{aligned}
    \d_t a^\e +\nabla \eik\cdot \nabla a^\e +
  \frac{1}{2}a^\e\Delta \eik &=i\frac{\e}{2}\Delta a^\e
  -i\e^{\kappa -1}f\(|a^\e|^2\)a^\e ,\\
a^\e_{\mid  t=0}&=a_0^\e. 
  \end{aligned}
\end{equation}
\begin{prop}\label{prop:asub}
  Let Assumptions~\ref{hyp:geom} and \ref{hyp:ana} be satisfied. Let
  $\kappa \ge 1$. For any $\e \in ]0,1]$, \eqref{eq:asub} has a
  unique solution $a^\e \in
  C^\infty([0,T]\times \R^n)\cap C([0,T];H^s)$ for any
  $s>n/2$. Moreover, $a^\e$ is bounded in 
  $L^\infty([0,T]; H^s)$ uniformly in $\e \in ]0,1]$, for any $s\ge
  0$. 
\end{prop}
\begin{proof}
  Using a mollification procedure, we see that it is enough to
  establish energy estimates for \eqref{eq:asub} in $H^s$, for any
  $s\ge 0$. Let $s>n/2$, and $\alpha\in \N^n$, with $s=|\alpha|$. Applying
  $\d^\alpha_x$ to \eqref{eq:asub}, we find:
\begin{equation}
  \label{eq:asubder}
    \d_t \d^\alpha_x a^\e +\nabla \eik\cdot \nabla \d^\alpha_x a^\e
    =i\frac{\e}{2}\Delta \d^\alpha_x a^\e 
-i\e^{\kappa -1}\d^\alpha_x\(f\(|a^\e|^2\)a^\e\) 
+ R_\alpha^\e, 
\end{equation}
where
\begin{equation*}
  R_\alpha^\e = \left[ \nabla \eik\cdot \nabla ,\d^\alpha_x \right]a^\e -
  \frac{1}{2}\d^\alpha_x \(a^\e\Delta \eik\). 
\end{equation*}
Take the inner product of \eqref{eq:asubder} with $\d^\alpha a^\e$,
and consider the real part:
the first term of the right hand side of \eqref{eq:asubder} vanishes,
and we have:
\begin{equation*}
  \begin{aligned}
    \frac{1}{2}\frac{d}{dt}\|\d^\alpha_x a^\e\|_{L^2}^2+ {\rm Re}\int_{\R^n}
  \overline{\d^\alpha_x a^\e}\nabla 
  \eik\cdot \nabla \d^\alpha_x a^\e \le &\e^{\kappa -1}
  \left\|f\(|a^\e|^2\)a^\e\right\|_{H^s} \left\|a^\e\right\|_{H^s}\\
  &   +\left\|R_\alpha^\e\right\|_{L^2}\left\|a^\e\right\|_{H^s} .
\end{aligned}
\end{equation*}
Notice that we have
\begin{align*}
  \left| {\rm Re}\int_{\R^n}
  \overline{\d^\alpha_x a^\e}\nabla 
  \eik\cdot \nabla \d^\alpha_x a^\e \right| &= \frac{1}{2}\left| \int_{\R^n}
  \nabla 
  \eik\cdot \nabla \left|\d^\alpha_x a^\e \right|^2\right|\\
& = \frac{1}{2}\left| \int_{\R^n}
  \left|\d^\alpha_x a^\e \right|^2 \Delta \eik  \right|\le C
\left\|a^\e\right\|_{H^s}^2, 
\end{align*}
since $\Delta \eik \in L^\infty([0,T]\times\R^n)$ from
Lemma~\ref{lem:hj}. Moser's inequality yields:
\begin{equation*}
  \left\|f\(|a^\e|^2\)a^\e\right\|_{H^s} \le C \(
  \left\|a^\e\right\|_{L^\infty}\) \left\|a^\e\right\|_{H^s}.
\end{equation*}
Summing over $\alpha$ such that $|\alpha|=s$, we infer:
\begin{equation*}
  \frac{d}{dt}\| a^\e\|_{H^s} \le C \(
  \left\|a^\e\right\|_{L^\infty}\) \left\|a^\e\right\|_{H^s} +
  \left\|R_\alpha^\e\right\|_{H^s}. 
\end{equation*}
Note that the above locally bounded map $C(\cdot)$ is independent of
$\e$ if and only if $\kappa \ge 1$. 
To apply Gronwall lemma, we need to estimate the last term: we use the
fact that the derivatives of order at least two of $\eik$ are bounded,
from Lemma~\ref{lem:hj}, to have:
\begin{equation*}
  \left\|R_\alpha^\e\right\|_{L^2} \le C \left\|a^\e\right\|_{H^s}.
\end{equation*}
We can then conclude by a continuity argument and Gronwall lemma: 
\begin{equation*}
  \| a^\e\|_{L^\infty([0,T];H^s)}\le C\(s, \left\|a_0^\e\right\|_{H^s} \). 
\end{equation*}
For the mollification procedure, this yields boundedness in the high
norm. Contraction in the small norm then follows easily with classical
arguments (see e.g. \cite{AlinhacGerard,Majda}), completing the proof
of Proposition~\ref{prop:asub}.
\end{proof}
\begin{cor}\label{cor:simplifsub}
 Let Assumptions~\ref{hyp:geom} and \ref{hyp:ana} be satisfied. Let
  $\kappa \ge 1$. Then
  \begin{equation*}
  \left\| a^\e -   \widetilde a^\e
  \right\|_{L^\infty([0,T];H^s)}\to 0 \quad \text{as }\e \to 0,\quad
  \forall s\ge 0, 
  \end{equation*}
where $\widetilde a^\e$ solves:
\begin{equation}\label{eq:atildeeps}
  \d_t \widetilde a^\e +\nabla \eik\cdot \nabla \widetilde a^\e +
  \frac{1}{2}\widetilde a^\e\Delta \eik =
  -i\e^{\kappa -1}f\(|\widetilde a^\e|^2\)\widetilde a^\e \quad ; \quad
\widetilde a^\e_{\mid  t=0}=a_0.
\end{equation}
\end{cor}
\begin{proof}
  The proof of Proposition~\ref{prop:asub} shows that $\widetilde a^\e \in
  C^\infty([0,T]\times \R^n)$, and that $\widetilde a^\e$ is bounded in
  $L^\infty([0,T]; H^s)$ uniformly in $\e \in ]0,1]$, for any $s\ge
  0$. Let $w^\e = a^\e -\widetilde a^\e$: it solves
\begin{align*}
  \d_t w^\e +\nabla \eik\cdot \nabla w^\e +
  \frac{1}{2}w^\e\Delta \eik &=i\frac{\e}{2}\Delta a^\e
  -i\e^{\kappa -1}\(F(a^\e) -F(\widetilde a^\e)\),\\ 
w^\e_{\mid  t=0}&=a_0^\e -a_0,
\end{align*}
where we have denoted $F(z)=f(|z|^2)z$. Proceeding as in the proof of
Proposition~\ref{prop:asub}, we have the following energy estimate:
\begin{equation*}
  \frac{d}{dt}\| w^\e\|_{H^s} \le C  \left\|w^\e\right\|_{H^s} +
  \e \left\|\Delta a^\e\right\|_{H^s} + \left\|F(a^\e) -F(\widetilde
    a^\e)\right\|_{H^s}.  
\end{equation*}
Since $H^s$ is an algebra and $F$ is $C^1$, the Fundamental Theorem of
Calculus yields:
\begin{equation*}
  \left\|F(a^\e) -F(\widetilde
    a^\e)\right\|_{H^s}\le C\(\left\| a^\e\right\|_{H^s},
  \left\|\widetilde a^\e\right\|_{H^s}\) \| w^\e\|_{H^s}.
\end{equation*}
Now since $a^\e$ and $\widetilde a^\e$ are bounded in
  $L^\infty([0,T]; H^{s+2})$ uniformly in $\e \in ]0,1]$, we have an
  estimate of the form:
\begin{equation*}
  \frac{d}{dt}\| w^\e\|_{H^s} \le C(s)  \left\|w^\e\right\|_{H^s} +
  C(s) \e .  
\end{equation*}
We conclude by Gronwall lemma, since $\|a_0^\e -a_0\|_{H^s} \to 0$
from Assumption~\ref{hyp:ana}.
\end{proof}
\begin{rema}
  If we assume moreover that like in \eqref{eq:CIbkw}, 
  \begin{equation*}
    a_0^\e = a_0 + \O(\e) \quad \text{in }H^s, \ \forall s\ge 0,
  \end{equation*}
then the above estimate can be improved to: 
\begin{equation*}
  \left\| a^\e -   \widetilde a^\e
  \right\|_{L^\infty([0,T];H^s)}=\O(\e),\quad \forall s\ge 0.
  \end{equation*}
\end{rema}
We have reduced the study of the asymptotic behavior of $u^\e$ to the
understanding of $\widetilde a^\e$. To complete the proof of
Proposition~\ref{prop:sub}, we resume the framework of
Section~\ref{sec:HJ}. With $x(t,y)$ given by the Hamilton flow
\eqref{eq:hamilton}, introduce the Jacobi determinant
\begin{equation*}
  J_t(y) ={\rm det}\nabla_y x(t,y). 
\end{equation*}
Denote
\begin{equation*}
  A^\e(t,y) := \widetilde a^\e \(t, x(t,y) \)\sqrt{J_t(y)}.
\end{equation*}
We see that so long as $y\mapsto x(t,y)$ defines a global
diffeomorphism (which is guaranteed for $t\in [0,T]$ by construction),
\eqref{eq:atildeeps} is equivalent to:  
\begin{equation*}
  \d_t A^\e = -i\e^{\kappa -1}f\(J_t(y)^{-1}\left|A^\e \right|^2\)
  A^\e\quad ; \quad A^\e(0,y)=a_0(y). 
\end{equation*}
This ordinary differential equation along the rays of geometrical
optics can be solved explicitly: since $f$ is real-valued, we see that
$\d_t |A^\e|^2 =0$, hence
\begin{equation*}
  A^\e(t,y) = a_0(y) \exp\(-i\e^{\kappa -1} \int_0^t
  f\(J_s(y)^{-1}\left|a_0(y) \right|^2\)ds\). 
\end{equation*}
Back to the function $\widetilde a^\e$, Proposition~\ref{prop:sub}
follows, with:
\begin{equation}\label{eq:exprG}
  \begin{aligned}
    a(t,x) &= \frac{1}{\sqrt{J_t(y(t,x))}}a_0\(y(t,x)\),\\
  G(t,x) &= -\int_0^t
  f\(J_s(y(t,x))^{-1}\left|a_0(y(t,x)) \right|^2\)ds.
  \end{aligned}
\end{equation}
One may wonder if this approach could be extended to some values
$\kappa <1$. Seek a solution of \eqref{eq:asub}. To have a simple ansatz
as in Proposition~\ref{prop:sub}, we would like to remove the
Laplacian in the limit $\e \to 0$ in \eqref{eq:asub}, and obtain the
analogue of Corollary~\ref{cor:simplifsub}. Following the same lines
as above, we find:
\begin{equation*}
  \widetilde a^\e (t,x) = a(t,x) e^{i\e^{\kappa-1}G(t,x)},
\end{equation*}
which is exactly the first formula in Proposition~\ref{prop:sub}. Now
recall that in Proposition~\ref{prop:asub}, we prove that $a^\e$ is
bounded in $H^s$, uniformly for $\e \in ]0,1]$; this property is used
to approximate $a^\e$ by $\widetilde a^\e$. But when $\kappa <1$,
$\widetilde a^\e$ is no longer uniformly bounded in $H^s$, because
what was a phase modulation for $\kappa \ge 1$ is now a rapid
oscillation. 
\smallbreak

This remark in a particular case reveals a much more general
phenomenon. When studying geometric optics in a supercritical r\'egime
(when $\kappa <1$ in the present context), distinguishing phase and
amplitude becomes a much more delicate issue. Suppose for instance
that we seek $u^\e = a^\e e^{i\phi^\e/\e}$, where $a^\e$ and $\phi^\e$
have asymptotic expansions as $\e \to 0$. All the terms
for $\phi^\e$ which are not $o(\e)$ are relevant, since $\phi^\e$ is
divided by $\e$. To determine these terms, it is not sufficient to
determine the leading order amplitude $\lim a^\e$ in general: because of
supercritical interactions, initial perturbations of the amplitude
may develop  non-negligible phase terms. To illustrate this discussion
in the above case, let $\kappa <1$ and $\widetilde a^\e$ solve
\eqref{eq:atildeeps} with initial data
\begin{equation*}
  \widetilde a^\e_{\mid t=0} =a_0 +\e^\gamma a_1,\quad \gamma >0,
\end{equation*}
where $a_0$ and $a_1$ are smooth and independent of $\e$. Integrating
\eqref{eq:atildeeps}, we find
\begin{equation*}
 \widetilde a^\e (t,x) =  \frac{1}{\sqrt{J_t(y(t,x))}}\(
 a_0\(y(t,x)\)+ \e^\gamma a_1\(y(t,x)\)\) e^{i\e^{\kappa -1}G^\e(t,x)},
\end{equation*}
where 
\begin{equation*}
  G^\e (t,x) = -\int_0^t
  f\(J_s(y(t,x))^{-1}\left|a_0(y(t,x))+ \e^\gamma a_1\(y(t,x)\) \right|^2\)ds.
\end{equation*}
We have the identity $G^\e = G + \e^\gamma G_1 +\e^{2\gamma}G_2$, for
the same $G$ as before, and $G_1,G_2$ independent of $\e$, but
depending on $a_1$. We infer:
\begin{equation*}
 \widetilde a^\e (t,x) \Eq \e 0  \frac{1}{\sqrt{J_t(y(t,x))}}
 a_0\(y(t,x)\) e^{i\e^{\kappa -1}(G+\e^\gamma G_1 +\e^{2\gamma}G_2) (t,x)}.
\end{equation*}
In particular, if $\kappa +\gamma \le 1$, we see that $G_1$ has to be
incorporated to describe the leading order behavior of $\widetilde
a^\e$. Since the three requirements $\kappa <1$, $\gamma>0$ and
$\kappa +\gamma \le 1$ can be met, we see that the small initial
perturbation $\e^\gamma a_1$ may produce relevant phase
perturbations. This example explains why in
Assumption~\ref{hyp:surcrit}, we require the asymptotic behavior of
the initial amplitude up to order $o(\e)$ and not only $o(1)$ (take
$\kappa=0$ and $\gamma =1$). 
We refer to
\cite{CheverryBullSMF} for a 
general explanation of this phenomenon, and to the next three sections
as far as Schr\"odinger equations are concerned.


\section{A hyperbolic point of view}
\label{sec:grenier}

In this section, we study \eqref{eq:nlssemi} in the case
$\kappa =0$, with no potential: 
\begin{equation}\label{eq:hs}
  i\e \d_t u^\e +\frac{\e^2}{2}\Delta u^\e =  f\(|u^\e|^2\)
  u^\e\quad ; \quad u^\e_{\mid t=0} = a_0^\e(x)e^{i\varphi_0(x)/\e}\, . 
\end{equation}
We recall the method introduced by E.~Grenier \cite{Grenier98},
which is valid for smooth nonlinearities which are defocusing and
cubic at the origin. Throughout this section, we assume the
following:
\begin{hyp}[Study of \eqref{eq:hs}]\label{hyp:hyper}
  We have $f'>0$. In addition,
  $\varphi_0\in H^\infty$, and there exists $a_0,a_1 \in H^\infty$
  such that
  \begin{equation*}
    a_0^\e = a_0 +\e a_1 +o(\e)\quad \text{in }H^s, \ \forall s\ge 0.
  \end{equation*}
\end{hyp}
Note that this assumption is closely akin to
Assumption~\ref{hyp:surcrit}: nevertheless, we do not make any
assumption on the 
momenta of $a_0$ and $a_1$, and  the initial phase is
bounded. We will see in Section~\ref{sec:surcrit} how to weaken this
assumption. 
\subsection{Grenier's approach}
\label{sec:hyper}
The principle is somehow to perform the usual WKB analysis ``the other
way round''. First, write the \emph{exact} solution as 
\begin{equation}\label{eq:ecrexacte}
  u^\e(t,x)= a^\e(t,x)e^{i\Phi^\e(t,x)/\e}\, ,
\end{equation}
where $\Phi^\e$ is real-valued. Then show that the ``amplitude'' $a^\e$
and the ``phase'' $\Phi^\e$ have asymptotic expansions as $\e\to 0$:
\begin{equation*}
  a^\e \sim a +\e a^{(1)}+\e^2 a^{(2)}+\ldots \quad ;\quad 
\Phi^\e \sim \phi +\e \phi^{(1)}+\e^2 \phi^{(2)}+\ldots
\end{equation*}
Introducing two unknown functions to solve one
equation yields a degree of freedom. The historical approach
\cite[Chap.~III]{LandauQ} consisted
in writing 
\begin{align*}
    \d_t \Phi^\e +\frac{1}{2}\left|\nabla \Phi^\e\right|^2 + f\(
    |a^\e|^2\)= \e^2 \frac{\Delta a^\e}{2 a^\e}\quad &; \quad
    \Phi^\e\big|_{t=0}=\varphi_0\, ,\\
\d_t a^\e +\nabla \Phi^\e \cdot \nabla a^\e +\frac{1}{2}a^\e
\Delta \Phi^\e  = 0\quad  &;\quad
a^\e\big|_{t=0}= a^\e_0\, . 
\end{align*}
Of course, this choice is not adapted when the amplitude $a^\e$
vanishes (see \cite{PGX93}), so it must be left out when $a^\e_0 \in
L^2(\R^n)$ in general. 
The approach introduced by
E.~Grenier consists in imposing:
\begin{equation}\label{eq:systexact0}
  \begin{aligned}
    \d_t \Phi^\e +\frac{1}{2}\left|\nabla \Phi^\e\right|^2 + f\(
    |a^\e|^2\)= 0\quad &; \quad
    \Phi^\e\big|_{t=0}=\varphi_0\, ,\\
\d_t a^\e +\nabla \Phi^\e \cdot \nabla a^\e +\frac{1}{2}a^\e
\Delta \Phi^\e  = i\frac{\e}{2}\Delta a^\e\quad & ;\quad
a^\e\big|_{t=0}= a^\e_0\, . 
  \end{aligned}
\end{equation}
Before recalling the results of \cite{Grenier98}, observe that if
$a^\e$ and $\Phi^\e$ are bounded in some sufficiently small Sobolev
spaces uniformly in $\e$, passing to the limit formally in
\eqref{eq:systexact0} yields:
\begin{equation}\label{eq:systlim}
  \begin{aligned}
    \d_t \phi +\frac{1}{2}\left|\nabla \phi\right|^2 + f\(
    |a|^2\)= 0\quad &; \quad
    \phi\big|_{t=0}=\varphi_0\, ,\\
\d_t a +\nabla \phi \cdot \nabla a +\frac{1}{2}a
\Delta \phi  = 0\quad & ;\quad
a\big|_{t=0}= a_0\, . 
  \end{aligned}
\end{equation}
We see that when the nonlinearity is
exactly cubic ($f(y)\equiv y$), $(\rho,v):=\( |a|^2,\nabla
\phi\)$ solves the  compressible, isentropic Euler equation 
\begin{equation}
  \label{eq:euler}
  \begin{aligned}
    \d_t v +v\cdot \nabla v + \nabla \rho= 0\quad &; \quad
    v\big|_{t=0}=\nabla \varphi_0\, ,\\
\d_t \rho  +\nabla\cdot (\rho v)  = 0\quad & ;\quad
\rho\big|_{t=0}= |a_0|^2\, . 
  \end{aligned}
\end{equation}
From this point of view, the formulation
\eqref{eq:systlim} is closely akin to the change of unknown function $\rho
\to \sqrt\rho$ introduced in \cite{MUK86} (see also
\cite{JYC90}) to study \eqref{eq:euler} when the initial density is
compactly supported, a situation more or less similar to the present
one. Note however that here, $a$ is complex-valued in general.
\smallbreak

Introducing the ``velocity'' $v^\e = \nabla \Phi^\e$,
\eqref{eq:systexact0} yields
\begin{equation}\label{eq:systexact}
  \begin{aligned}
    \d_t v^\e +v^\e \cdot \nabla v^\e + 2 f'\(
    |a^\e|^2\) \operatorname{Re}\(\overline{a^\e}\nabla a^\e\)=
    0\quad &; \quad 
    v^\e\big|_{t=0}=\nabla \phi_0\, ,\\
\d_t a^\e +v^\e\cdot \nabla a^\e +\frac{1}{2}a^\e
\nabla\cdot v^\e  = i\frac{\e}{2}\Delta a^\e\quad & ;\quad
a^\e\big|_{t=0}= a^\e_0\, . 
  \end{aligned}
\end{equation}
Separate real and
imaginary parts of $a^\e$, $a^\e = 
a_1^\e + ia_2^\e$. Then we have 
\begin{equation}
  \label{eq:systhyp}
  \d_t \bu^\e +\sum_{j=1}^n A_j(\bu^\e)\d_j \bu^\e
  = \frac{\e}{2} L 
  \bu^\e\, , 
\end{equation}
\begin{equation*}
  \text{with}\quad \bu^\e = \left(
    \begin{array}[l]{c}
       a_1^\e \\
       a_2^\e \\
       v^\e_1 \\
      \vdots \\
       v^\e_n
    \end{array}
\right)\quad , \quad L = \left(
  \begin{array}[l]{ccccc}
   0  &-\Delta &0& \dots & 0   \\
   \Delta  & 0 &0& \dots & 0  \\
   0& 0 &&0_{n\times n}& \\
   \end{array}
\right),
\end{equation*}
\begin{equation*}
  \text{and}\quad A(\bu,\xi)=\sum_{j=1}^n A_j(\bu)\xi_j
= \left(
    \begin{array}[l]{ccc}
      v\cdot \xi & 0& \frac{a_1 }{2}\,^{t}\xi \\ 
     0 &  v\cdot \xi & \frac{a_2}{2}\,^{t}\xi \\ 
     2f'  a_1 \, \xi
     &2f'  a_2\, \xi &  v\cdot \xi I_n 
    \end{array}
\right),
\end{equation*}
where $f'$ stands for $f'(|a_1|^2+|a_2|^2)$. 
The matrix $A(\bu,\xi)$ can be symmetrized by 
\begin{equation}\label{eq:symetriseur}
  S=\left(
    \begin{array}[l]{cc}
     I_2 & 0\\
     0& \frac{1}{4f'}I_n
    \end{array}
\right),
\end{equation}
which is symmetric and positive since $f'>0$. For an integer
$s>2+n/2$, we bound $(S 
\d_{x}^{\alpha } \bu^\e 
, \d_{x}^{\alpha }\bu^\e )$ 
where $\alpha $ is a  multi index of length $\le s$, and
$(\cdot ,\cdot)$ is the usual  $L^{2}$ scalar product. We have
\begin{equation*}
\frac{d}{dt}\(S \d_{x}^{\alpha } \bu^\e ,
  \d_{x}^{\alpha } \bu^\e\) 
= \(\d_{t} S  \d_{x}^{\alpha } \bu^\e , \partial
_{x}^{\alpha } \bu^\e\) 
  + 2 \( S \d_{t}  \d_{x}^{\alpha } \bu^\e , \partial
  _{x}^{\alpha } \bu^\e\) 
 \end{equation*}
since $S$ is symmetric. For the first term, we consider the lower
$n\times n$ block:
\begin{equation*}
\(\d_{t} S  \d_{x}^{\alpha } \bu^\e ,
  \d_{x}^{\alpha } \bu^\e\) 
\le  \left\|\frac{1}{f'}\d_t\(f'\( | a_1^\e|^2 + 
| a_2^\e|^2\)\)\right\|_{L^\infty}\(S  \d_{x}^{\alpha } \bu^\e ,
  \d_{x}^{\alpha } \bu^\e\)\, .
\end{equation*}
So long as $\|\bu^\e\|_{L^\infty}\le 2\|a_0^\e\|_{L^\infty}$, we have:
\begin{equation*}
  f'\( |a_1^\e|^2 + 
| a_2^\e|^2\) \ge \inf\left\{ f'(y)\ ;\ 0\le y\le
  4\sup_{0<\e\le 1}\|a_0^\e\|_{L^\infty}^2\right\} =\delta_n>0\, ,
\end{equation*}
where $\delta_n$ is now fixed, since $f'$ is continuous with
$f'>0$. We infer, 
\begin{equation*}
 \left\|\frac{1}{f'}\d_t\(f'\( |a_1^\e|^2 + 
| a_2^\e|^2\)\)\right\|_{L^\infty}\le C \|\bu^\e\|_{H^s}\, ,
\end{equation*}
where we used Sobolev embeddings and
\eqref{eq:systhyp}. 
For the second term we use 
\begin{equation*} 
\( S \d_{t}  \d_{x}^{\alpha } \bu^\e , \d_{x}^{\alpha } \bu^\e\)
= \frac{\e}{2}\(S L( \d_{x}^{\alpha } \bu^\e )  , \d_{x}^{\alpha }
\bu^\e\)
  -  
\Big( S  \d_{x}^{\alpha } \Big(\sum _{j=1}^{n} A_j(\bu^\e) \d_{j}
  \bu^\e  \Big) , 
 \d_{x}^{\alpha } \bu^\e\Big).
\end{equation*} 
We notice that $SL$
is a skew-symmetric second order operator, so the first term is zero. 
For the second term, use the symmetry of $SA_j(\bu^\e)$  and usual
estimates on commutators  to get finally:
\begin{equation*}
\frac{d}{dt} \sum _{|\alpha | \le s} \(S \d_{x}^{\alpha } \bu^\e,
\d_{x}^{\alpha } \bu^\e\) \le C\(\left\|\bu^\e\right\|_{H^s}\)
\sum _{|\alpha | \le s} \(S \d_{x}^{\alpha } \bu^\e,
\d_{x}^{\alpha } \bu^\e\)\, ,
\end{equation*}
for $s > 2+d/2 $.
Gronwall lemma along with a
continuity argument yield: 
\begin{prop}[\cite{Grenier98}, Th.~1.1]\label{prop:p}
Let Assumption~\ref{hyp:hyper} be satisfied. Let
$s>2+n/2$. There exist $T_s>0$ independent of $\e\in ]0,1]$  and
$u^\e = a^\e e^{i\Phi^\e/\e}$ solution to
\eqref{eq:hs} on 
$[0,T_s]$. Moreover, $a^\e$ and
$\Phi^\e$ are bounded 
in $L^\infty([0,T_s];H^s)$,
uniformly in $\e\in ]0,1]$. 
\end{prop}
The solution to \eqref{eq:systexact0} formally converges to the
solution  of \eqref{eq:systlim}. Under Assumption~\ref{hyp:hyper}, 
\eqref{eq:systlim} has a unique solution $(a,\phi)\in
L^\infty([0,T_*];H^m)^2$ for any $m>0$ for some $T_*>0$ 
independent of $m$ (see e.g. \cite{AlinhacGerard,Majda}). We
infer: 
\begin{prop}\label{prop:estprec}
Let $s\in \N$. Then $T_s\ge T_*$, and there exists $C_s$ independent of
$\e$ such that for 
every $0\le t\le T_*$, 
\begin{equation}\label{eq:tpetit}
  \| a^\e (t)- a(t)\|_{H^s}\le C_s \e\quad ;\quad
  \| \Phi^\e(t) - \phi(t)\|_{H^s} \le C_s \e t . 
\end{equation}
\end{prop}
\begin{proof}
  We keep the same notations as above,
  \eqref{eq:systhyp}. Denote by $\bv$
  the analog of $\bu^\e$ 
  corresponding to $(a,\phi)$. We have
  \begin{equation*}
\d_t \(\bu^\e-\bv\) +\sum_{j=1}^n A_j(\bu^\e)\d_j \(\bu^\e-\bv\) +
\sum_{j=1}^n \(A_j(\bu^\e)-A_j(\bv)\)\d_j \bv  
  = \frac{\e}{2} L 
  \bu^\e    \, .
  \end{equation*}
Keeping the symmetrizer $S$ corresponding to $\bu^\e$, we can do
similar computations to  the previous ones. Note that we
know that $\bu^\e$ and $\bv$ are bounded in
$L^\infty([0,\min(T_s,T_*)];H^s)$. Denote ${\bw}^\e 
=\bu^\e-\bv$. Writing $L 
  \bu^\e = L {\bw}^\e + L\bv$, the term $L {\bw}^\e$ disappears from
  the energy estimate, and we get, for $s>2+n/2$: 
\begin{align*}
\frac{d}{dt} \sum _{|\alpha | \le s} \(S \d_{x}^{\alpha } \bw^\e,
\d_{x}^{\alpha } \bw^\e\) \le & \, C\( \|\bu^\e\|_{H^s},
\|\bv\|_{H^{s+2}}\) 
\sum _{|\alpha | \le s} \(S \d_{x}^{\alpha } \bw^\e,
\d_{x}^{\alpha } \bw^\e\)\\
&+ \e \|\bv\|_{H^{s+2}} \| \bw^\e(t)\|_{H^s}\, .
\end{align*}
Gronwall lemma and a continuity argument show that
$\bw^\e$ (hence $\bu^\e$) is defined on $[0,T_*]$. By
Assumption~\ref{hyp:hyper}, $\|\bw^\e_{\mid t=0}\|_{H^s} = \O(\e)$,
and we get:
\begin{equation*}
  \left\| \bw^\e\right\|_{L^\infty([0,T_*];H^s)}=\O(\e ). 
\end{equation*}
The estimate for the phase (and not only its gradient) then follows
from the above estimate and the integration in time of
\eqref{eq:systexact0}-\eqref{eq:systlim}.  
\end{proof}
Proposition~\ref{prop:estprec} yields an approximation of $u^\e$
\emph{for small times only}:
\begin{align*}
  \big\| u^\e(t) - & a(t)e^{i\phi(t)/\e}\big\|_{L^2}=\left\|
  a^\e(t)e^{i\Phi^\e(t)/\e}- a(t)e^{i\phi(t)/\e}\right\|_{L^2}\\
&=\O\(\left\|
  a^\e(t)- a(t)\right\|_{L^2}+\left\|
  e^{i\Phi^\e(t)/\e}- e^{i\phi(t)/\e}\right\|_{L^\infty }\left\|
  a(t)\right\|_{L^2} \)\\
&  = \O(\e)+\O(t). 
\end{align*}
For times of order $\O(1)$, the corrector $a_1$ must be taken into
account: 
\begin{prop}\label{prop:correc}
  Let Assumption~\ref{hyp:hyper} be satisfied. Define
  $(a^{(1)},\phi^{(1)})$ by 
\begin{equation*}
  \begin{aligned}
    \d_t \phi^{(1)} +\nabla \phi \cdot \nabla \phi^{(1)} +
    2\operatorname{Re}\(\overline a a^{(1)}\)f'\( |a|^2\)&=0,\\
   \d_t a^{(1)} +\nabla\phi\cdot \nabla a^{(1)} + \nabla
   \phi^{(1)}\cdot \nabla a + \frac{1}{2} a^{(1)}\Delta \phi
   +\frac{1}{2}a\Delta \phi^{(1)}      &= \frac{i}{2}\Delta a,\\
\phi^{(1)}\big|_{t=0}=0\quad ; \quad a^{(1)}\big|_{t=0}=a_1.
  \end{aligned}
\end{equation*}
Then $a^{(1)},\phi^{(1)}\in
L^\infty([0,T_*];H^s)$ for every $s\ge 0$, and
\begin{equation*}
  \|a^\e - a - \e a^{(1)}\|_{L^\infty([0,T_*];H^s)}+ \|\Phi^\e - \phi - \e
  \phi^{(1)}\|_{L^\infty([0,T_*];H^s)} \le C_s\e^2,\quad \forall s\ge 0\, .
\end{equation*}
\end{prop}
The proof is a straightforward consequence of the above analysis, and
is given in \cite{Grenier98}. Despite 
the notations, it seems unadapted to consider $\phi^{(1)}$ as being
part of the phase. Indeed, we infer from Proposition~\ref{prop:correc}
that 
\begin{equation*}
  \left\|u^\e - a e^{i\phi^{(1)}}
    e^{i\phi/\e}\right\|_{L^\infty([0,T_*];L^2\cap L^\infty)}= \O(\e). 
\end{equation*}
Relating this information to the WKB methods presented in
the introduction, we would have:
\begin{equation*}
  \a_0 = a e^{i\phi^{(1)}}.
\end{equation*}
Since $\phi^{(1)}$ depends on $a_1$ while $a$ does not, we retrieve
the fact that in super-critical r\'egimes, the leading order amplitude
in WKB methods depends on the initial first corrector $a_1$. 
\begin{rema}
  The term $e^{i\phi^{(1)}}$ does not appear in the Wigner measure
  of $a e^{i\phi^{(1)}} e^{i\phi/\e}$. Thus, from the point of view of
    Wigner measures, the asymptotic behavior of the exact solution is
    described by the Euler-type system \eqref{eq:systlim}. 
\end{rema}

\begin{rema}[Introducing an isotropic harmonic potential]
The above method makes it possible to consider the semi-classical of
\eqref{eq:nlssemi} when $V(t,x)=\frac{1}{2}|x|^2$ is an isotropic harmonic
potential, and Assumption~\ref{hyp:hyper} is satisfied. Let 
\begin{equation*}
U^\e(t,x) =
\frac{1}{(1+t^2)^{n/4}}e^{i\frac{t}{1+t^2}\frac{|x|^2}{2\e}}u^\e \(
\arctan t , \frac{x}{\sqrt{1+t^2}}\)\, .
\end{equation*}
Then $U^\e$ solves:
\begin{equation*}\left\{
\begin{aligned}
i\e \d_t U^\e + \frac{\e^2}{2}\Delta U^\e& = 
\frac{1}{1+t^2}f\(\(1+t^2\)^{n/2}|U^\e|^2\)U^\e\, , \\
U^\e(0,x)&= a_0^\e(x)e^{i\varphi_0(x)/\e} .
\end{aligned}\right.
\end{equation*}
We can then proceed as above. The only difference is the presence of
time in the nonlinearity, which changes very little the analysis. 
\end{rema}

\begin{rema}[Momenta]\label{rema:moment}
If in Assumption~\ref{hyp:hyper}, we replace $H^s$ with 
\begin{equation*}
  \Sigma^s =H^s \cap \F (H^s) =\left\{ w \in L^2\ ; \
    (1-\Delta)^{k/2}\< x\>^{s-k}w\in L^2, \   0\le k\le s\right\},
\end{equation*}
then the above analysis can be repeated in $\Sigma^s$. The main
difference is due to the commutations of the
powers of $x$ with the differential operators; it is easy to check
that they introduce semilinear terms, which can be treated as
source terms when applying  Gronwall lemma. 
\end{rema}

\subsection{Remarks about some conserved quantities}
\label{sec:conserved}
Consider the case of the cubic, defocusing Schr\"odinger equation:
$f(y)\equiv y$. 
Recall three important evolution laws for \eqref{eq:nlssemi}:
\begin{align*}
  \text{Mass: }& \frac{d}{dt}\|u^\e(t)\|_{L^2}=0\, .\\
\text{Energy: }& \frac{d}{dt}\(\|\e\nabla_x u^\e\|_{L^2}^2
+ \|u^\e\|_{L^4}^4 \)=0\, .\\
\text{Momentum: }&\frac{d}{dt}\operatorname{Im}\int \overline {u^\e}(t,x)
\e \nabla_x u^\e(t,x)dx =0\, .\\ 
\text{Pseudo-conformal law: }& \frac{d}{dt}\(\|J^\e(t)
u^\e\|_{L^2}^2 
+ t^2\|u^\e\|_{L^{4}}^{4}
\)=t (2-n)\|u^\e\|_{L^{4}}^{4} ,
\end{align*}
where $J^\e(t) = x + i\e t\nabla_x$. These evolutions are deduced
from the usual ones ($\e =1$, see e.g. \cite{CazCourant,Sulem}) via
the scaling 
$\psi(t,x)= u(\e t,\e x)$. Using \eqref{eq:ecrexacte} and passing to
the limit formally in the above formulae yields:
\begin{align*}
   \frac{d}{dt}\|a(t)\|_{L^2}&=0\, .\\
 \frac{d}{dt}\int \( |a(t,x)|^2 |\nabla \phi(t,x)|^2 +
|a(t,x)|^4\)dx &=0\, .\\
\frac{d}{dt}\int |a(t,x)|^2 \nabla \phi(t,x)dx &=0\, .\\
 \frac{d}{dt}\int \( \left| \(x -t \nabla
  \phi(t,x)\)a(t,x)\right|^2 +t^2 |a(t,x)|^4 
\)&dx=\\
=(2-n)t &\int |a(t,x)|^4dx\, .
\end{align*}
Note that we also
have the conservation (\cite{CN}):
\begin{equation*}
   \frac{d}{dt}\operatorname{Re}\int \overline{u^\e}(t,x)
 J^\e(t)u^\e(t,x) dx=0 \, ,
\end{equation*}
which yields:
\begin{equation*}
  \frac{d}{dt}\int \(\(x -t \nabla
  \phi(t,x)\)|a(t,x)|^2 \)dx =0\, . 
\end{equation*}
All these expressions involve only $(|a|^2, \nabla \phi)$, that
is, the solution of \eqref{eq:euler}. We thus
retrieve formally some evolution laws
for the Euler equation.


\section{Introducing subquadratic potential and initial phase}
\label{sec:surcrit}
In this section, we prove Theorem~\ref{theo:BKWV}. 
First, we point out that the uniqueness for $u^\e$ in $C([0,T_*];H^s)$
is straightforward for $s>n/2$. We thus have to prove that there
exists such a solution, and that it is smooth. 
\smallbreak

As suggested by
the statements of Theorem~\ref{theo:BKWV}, the idea consists in
resuming Grenier's method, and in writing the phase $\Phi^\e$ as 
\begin{equation*}
  \Phi^\e = \eik +\phi^\e.
\end{equation*}
We take $\phi^\e$ as a new unknown function. Recall that the
system~\eqref{eq:systexact0} reads, with the present notations:
\begin{align*}
    \d_t \Phi^\e +\frac{1}{2}\left|\nabla \Phi^\e\right|^2 + V+ f\(
    |a^\e|^2\)= 0\quad &; \quad
    \Phi^\e\big|_{t=0}=\phi_0\, ,\\
\d_t a^\e +\nabla \Phi^\e \cdot \nabla a^\e +\frac{1}{2}a^\e
\Delta \Phi^\e  = i\frac{\e}{2}\Delta a^\e\quad & ;\quad
a^\e\big|_{t=0}= a^\e_0\, . 
\end{align*}
This system becomes, in terms of $\phi^\e$, and given \eqref{eq:eik}: 
\begin{equation}\label{eq:modif0}
\begin{aligned}
    \d_t \phi^\e +\frac{1}{2}\left|\nabla \phi^\e\right|^2 +\nabla
    \eik\cdot \nabla \phi^\e+ f\(
    |a^\e|^2\)&= 0,\\
    \d_t a^\e +\nabla \phi^\e \cdot \nabla a^\e +\nabla \eik \cdot \nabla
    a^\e +\frac{1}{2}a^\e 
\Delta \phi^\e +\frac{1}{2}a^\e 
\Delta \eik &= i\frac{\e}{2}\Delta a^\e,\\
\phi^\e\big|_{t=0}=0\quad ;\quad a^\e\big|_{t=0}&= a^\e_0\, .
\end{aligned} 
\end{equation}
Like in Section~\ref{sec:hyper}, we 
work with $v^\e =\nabla \phi^\e$ instead of $\phi^\e$, to begin
with. The new terms are the factors where $\eik$ is present. The point
is to check that they are semilinear perturbations, which can be treated as
source terms in view of Gronwall lemma. Again, separate real and
imaginary parts of $a^\e$, $a^\e = 
a_1^\e + ia_2^\e$, and introduce:
\begin{equation*}
  \bu^\e = \left(
    \begin{array}[l]{c}
       a_1^\e \\
       a_2^\e \\
       v^\e_1 \\
      \vdots \\
       v^\e_n
    \end{array}
\right)\quad , \quad L = \left(
  \begin{array}[l]{ccccc}
   0  &-\Delta &0& \dots & 0   \\
   \Delta  & 0 &0& \dots & 0  \\
   0& 0 &&0_{n\times n}& \\
   \end{array}
\right),
\end{equation*}
\begin{equation*}
  \text{and}\quad A(\bu,\xi)=\sum_{j=1}^n A_j(\bu)\xi_j
= \left(
    \begin{array}[l]{ccc}
      v\cdot \xi & 0& \frac{a_1 }{2}\,^{t}\xi \\ 
     0 &  v\cdot \xi & \frac{a_2}{2}\,^{t}\xi \\ 
     2f'  a_1 \, \xi
     &2f'  a_2\, \xi &  v\cdot \xi I_n 
    \end{array}
\right),
\end{equation*}
where $f'$ stands for $f'(|a_1|^2+|a_2|^2)$. Instead of
\eqref{eq:systhyp}, we now have a system of the form
\begin{equation}
  \label{eq:systhypV}
  \d_t \bu^\e +\sum_{j=1}^n A_j(\bu^\e)\d_j \bu^\e
  +\sum_{j=1}^n B_j(\nabla \eik)\d_j \bu^\e + M\(\nabla^2
  \eik\)\bu^\e= \frac{\e}{2} L  
  \bu^\e\, , 
\end{equation}
where the matrices $B_j$ depend linearly on their argument, and the matrix
$M$ is smooth, locally bounded. The quasilinear part of
\eqref{eq:systhypV} is the same as in
Section~\ref{sec:hyper}, and involves the matrices $A_j$. In
particular, we keep the same symmetrizer $S$ given by
\eqref{eq:symetriseur}. The matrices
$B_j$ have a semilinear contribution, as we see below. The term
corresponding to the matrix $M$ can obviously be considered as a
source term, since $\eik$ is subquadratic. 
\smallbreak

Let $s$ be an integer, 
$s>2+n/2$, and let $\alpha $ be a  multi index of length $\le s$. We
have: 
\begin{equation*}
\frac{d}{dt}\(S \d_{x}^{\alpha } \bu^\e ,
  \d_{x}^{\alpha } \bu^\e\) 
= \(\d_{t} S  \d_{x}^{\alpha } \bu^\e , \partial
_{x}^{\alpha } \bu^\e\) 
  + 2 \( S \d_{t}  \d_{x}^{\alpha } \bu^\e , \partial
  _{x}^{\alpha } \bu^\e\) 
 \end{equation*}
since $S$ is symmetric. For the first term, we consider the lower
$n\times n$ block:
\begin{equation}\label{eq:estdtS}
\(\d_{t} S  \d_{x}^{\alpha } \bu^\e ,
  \d_{x}^{\alpha } \bu^\e\) 
\le  \left\|\frac{1}{f'}\d_t\(f'\( | a_1^\e|^2 + 
| a_2^\e|^2\)\)\right\|_{L^\infty}\(S  \d_{x}^{\alpha } \bu^\e ,
  \d_{x}^{\alpha } \bu^\e\)\, .
\end{equation}
We consider times not larger than $T$ given by Lemma~\ref{lem:hj}, so
that the function $\eik$ remains smooth and subquadratic. 
So long as $\|\bu^\e\|_{L^\infty}\le 2\|a_0^\e\|_{L^\infty}$, we have:
\begin{equation*}
  f'\( |a_1^\e|^2 + 
| a_2^\e|^2\) \ge \inf\left\{ f'(y)\ ;\ 0\le y\le
  4\sup_{0<\e\le 1}\|a_0^\e\|_{L^\infty}^2\right\} =\delta_n>0\, .
\end{equation*}
We infer, 
\begin{equation}\label{eq:dtS}
 \left\|\frac{1}{f'}\d_t\(f'\( |a_1^\e|^2 + 
| a_2^\e|^2\)\)\right\|_{L^\infty}\le C \(\|\bu^\e\|_{H^s}+
\|x\bu^\e\|_{H^{s-1}}\)\, , 
\end{equation}
for some locally bounded map $C(\cdot)$. We used Sobolev embeddings,
\eqref{eq:systhypV} and Lemma~\ref{lem:hj}: the terms $B_j$ are
sublinear in $x$, hence the norm $\|x\bu^\e\|_{H^{s-1}}$ which we did
not consider in Section~\ref{sec:hyper}. We emphasize that this
estimate explains why we assume $s>2+n/2$, and not only $s>1+n/2$:
we control $\d_t \bu^\e$ in $L^\infty$ using  \eqref{eq:systhypV}, 
so we need to estimate $L \bu^\e$ in $L^\infty$. For all the other
terms, $s>1+n/2$ would suffice. This also explains why we wrote
$\|x\bu^\e\|_{H^{s-1}}$ and not $\|x\bu^\e\|_{H^{s}}$.  
For the second term we use 
\begin{align*} 
\( S \d_{t}  \d_{x}^{\alpha } \bu^\e , \d_{x}^{\alpha } \bu^\e\)
=& \frac{\e}{2}\(S L( \d_{x}^{\alpha } \bu^\e )  , \d_{x}^{\alpha }
\bu^\e\)
  -  
\Big( S  \d_{x}^{\alpha } \Big(\sum _{j=1}^{n} A_j(\bu^\e) \d_{j}
  \bu^\e  \Big) , 
 \d_{x}^{\alpha } \bu^\e\Big)\\
- \Big( S  \d_{x}^{\alpha } \Big(\sum _{j=1}^{n} B_j&(\nabla\eik) \d_{j}
  \bu^\e  \Big) , 
 \d_{x}^{\alpha } \bu^\e\Big)- \Big( S  \d_{x}^{\alpha }
 \Big(M(\nabla^2 \eik) 
  \bu^\e  \Big) , 
 \d_{x}^{\alpha } \bu^\e\Big).
\end{align*}   
The first two terms of the right hand side are handled in the same way
as in Section~\ref{sec:hyper}: the first one is zero, and the second
can be estimated by: 
\begin{equation}\label{eq:estQL}
 \Big( S  \d_{x}^{\alpha } \Big(\sum _{j=1}^{n} A_j(\bu^\e) \d_{j}
  \bu^\e  \Big) , 
 \d_{x}^{\alpha } \bu^\e\Big)\le C\(\left\|\bu^\e\right\|_{H^s}\)
\sum _{|\alpha | \le s} \(S \d_{x}^{\alpha } \bu^\e,
\d_{x}^{\alpha } \bu^\e\),
\end{equation}
where we keep the convention that $C(\cdot)$ is a locally bounded
map. Let us briefly explain this quasilinear estimate. First, write
\begin{equation*} 
\begin{aligned}
\( S  \partial _{x}^{\alpha } \( A_j(\bu^\e) \d_j
 \bu^\e  \) ,   
 \partial _{x}^{\alpha } \bu^\e\)&
= 
\(S   A_j(\bu^\e) \d_j \partial _{x}^{\alpha
 } \bu^\e, 
 \partial _{x}^{\alpha } \bu^\e\) \\
 +& \(S \(  \partial _{x}^{\alpha } (
 A_j(\bu^\e) \d_j \bu^\e  ) - 
  A_j(\bu^\e) \d_j \partial _{x}^{\alpha }
 \bu^\e  \), 
 \partial _{x}^{\alpha } \bu^\e\).
\end{aligned}
\end{equation*}
By symmetry of $SA_j(\bu^\e)$, 
\begin{equation*}
\begin{aligned}
\(S  A_j(\bu^\e) \d_j  \partial_{x}^{\alpha } \bu^\e, 
 \partial _{x}^{\alpha } \bu^\e\)
 = & -  \(\d_j (S A_j(\bu^\e))  \d_{x}^{\alpha } \bu^\e,
 \d_{x}^{\alpha } \bu^\e\)\\
 &-  \(S A_j(\bu^\e) \d_j  \partial _{x}^{\alpha } \bu^\e,
 \partial _{x}^{\alpha } \bu^\e\) .
\end{aligned} 
\end{equation*} 
We infer: 
\begin{align*} 
\left| \(S   A_j(\bu^\e) \d_j \partial _{x}^{\alpha } \bu^\e,
 \d_{x}^{\alpha } \bu^\e\) \right|&
\le
 \left\| \d_j\(S A_j (\bu^\e)\)
 \right\|_{L^{\infty }} 
\left\|\d_{x}^{\alpha } \bu^\e \right\|_{L^{2}}^{2}\\
&\le 
C\( \left\| \bu^\e \right\|_{L^{\infty }}\) \left\| \nabla _{x} \bu^\e
 \right\|_{L^{\infty }} 
\left\|\d_{x}^{\alpha } \bu^\e \right\|_{L^{2}}^{2} . 
\end{align*}
The usual estimates on commutators (see e.g. \cite{Majda}) lead to
\begin{equation*}
\left|\(S  \(  \partial_{x}^{\alpha } \(
  A_j(\bu^\e) \d_j \bu^\e  \) - 
  A_j(\bu^\e) \d_j \partial_{x}^{\alpha } \bu^\e \),
 \partial _{x}^{\alpha } \bu^\e\) \right|
\le C\(\left\|\bu^\e\right\|_{H^s}\) \left\|\bu^\e\right\|_{H^s}^{2},
\end{equation*}
and \eqref{eq:estQL} follows, since we consider times where 
 $S^{-1}$ is bounded. 
\smallbreak

For the third term of $\( S \d_{t}  \d_{x}^{\alpha } \bu^\e ,
\d_{x}^{\alpha } \bu^\e\)$, write: 
\begin{align*}
 \( S  \d_{x}^{\alpha } \(B_j(\nabla\eik) \d_{j}
  \bu^\e  \) , 
 \d_{x}^{\alpha } \bu^\e\) =& \int  S  B_j(\nabla\eik) \d_{j}\d_{x}^{\alpha } 
  \bu^\e\d_{x}^{\alpha } \bu^\e dx\\
&+\int  S  \left[\d_{x}^{\alpha } ,B_j(\nabla\eik) \d_{j}\right]
  \bu^\e\d_{x}^{\alpha } \bu^\e dx.  
\end{align*}
For the first term of the right hand side, an integration by parts
yields:
\begin{align}
 \left|\int  S  B_j(\nabla\eik) \d_{j}\d_{x}^{\alpha } 
  \bu^\e\d_{x}^{\alpha } \bu^\e dx\right| &\le \left\|
  \d_j\(  S  B_j(\nabla\eik)
  \)\right\|_{L^\infty}\|\bu^\e\|_{H^s}^2 \label{eq:16h47} \\ 
&\le C(\|a^\e \|_{L^\infty})\left\|\<x\> a^\e \nabla a^\e
  \right\|_{L^\infty}\|\bu^\e\|_{H^s}^2\notag \\ 
\le C(\|\bu^\e\|_{L^\infty})
\Big(\left\|\bu^\e\right\|_{L^\infty}&+\left\|x\bu^\e\right\|_{L^\infty} 
  + \left\|\nabla\bu^\e\right\|_{L^\infty}\Big)^2\|\bu^\e\|_{H^s}^2,\notag
\end{align}
where we have used Lemma~\ref{lem:hj}.  Again
from Lemma~\ref{lem:hj}, the commutator 
\begin{equation*}
  \left[\d_{x}^{\alpha } ,B_j(\nabla\eik) \d_{j}\right]
\end{equation*}
is a differential operator of degree $\le s$, with bounded
coefficients. We infer:
\begin{align*}
  \Big| \Big( S  \d_{x}^{\alpha } \Big(\sum _{j=1}^{n} B_j(\nabla\eik) \d_{j}
  \bu^\e  \Big) , 
 \d_{x}^{\alpha } \bu^\e\Big)\Big| 
\le C(\|\bu^\e \|_{H^s} +\|x\bu^\e \|_{H^{s-1}} )
  \|\bu^\e\|_{H^s}^2.
\end{align*}
We have obviously
\begin{align*}
  \Big| \Big( S  \d_{x}^{\alpha } \Big(M(\nabla^2\eik) 
  \bu^\e  \Big) , 
 \d_{x}^{\alpha } \bu^\e\Big)\Big| \le C(\|\bu^\e \|_{L^\infty})
  \|\bu^\e\|_{H^s}^2.
\end{align*}
This yields:
\begin{equation}\label{eq:10h45}
  \frac{d}{dt}\(S \d_{x}^{\alpha } \bu^\e ,
  \d_{x}^{\alpha } \bu^\e\)  \le C
  \(\left\|\bu^\e\right\|_{H^s}+\left\|x\bu^\e\right\|_{H^{s-1}} 
  \)\|\bu^\e\|_{H^s}^2, 
\end{equation}
where the map $C(\cdot)$ is locally bounded. 
We now have to
bound $x\bu^\e$ in $H^{s-1}$  to close our family of
estimates: we consider
\begin{equation*}
\frac{d}{dt}\(S \d_{x}^{\beta } (x_k\bu^\e) ,
  \d_{x}^{\beta } (x_k\bu^\e)\);\quad 1\le k\le n, \ |\beta|\le s-1.
 \end{equation*}
We can proceed as above, replacing $\bu^\e$ with $x_k\bu^\e$:
$x_k\bu^\e$ solves almost the same equation as $\bu^\e$, and we must
control some commutators. 
\begin{align*}
 \d_t (x_k\bu^\e) +\sum_{j=1}^n A_j(\bu^\e)\d_j &(x_k \bu^\e)
  +\sum_{j=1}^n B_j(\nabla \eik)\d_j (x_k \bu^\e)+ \\
+ M\(\nabla^2
  \eik\)x_k\bu^\e
= \frac{\e}{2} L  
  (x_k \bu^\e)
&+ A_k (\bu^\e) \bu^\e + B_k(\nabla \eik)\bu^\e +
  \frac{\e}{2}[x_k,L]\bu^\e .
\end{align*}
The term $A_k (\bu^\e) \bu^\e$ is harmless. The term $B_k(\nabla
\eik)\bu^\e$ is controlled by $\<x\> \bu^\e$, since $\eik$ is
subquadratic: this is a (semi)linear perturbation. Finally,
\begin{equation*}
 [x_k,L] = \left(
  \begin{array}[l]{ccccc}
   0  &2\d_k &0& \dots & 0   \\
   -2\d_k  & 0 &0& \dots & 0  \\
   0& 0 &&0_{n\times n}& \\
   \end{array}
\right).
\end{equation*}
Now we only have to notice that estimating $x_k \bu^\e$ does not
involve extra regularity or extra momenta. In the above computations,
the first time we needed to consider momenta was for \eqref{eq:dtS}:
we need exactly the same estimate now, since it is due to the
symmetrizer, which remains the same. The same remark is valid for
\eqref{eq:16h47}. For $\beta$ a multi index of length $\le s-1$, we
find:
\begin{equation}\label{eq:10h46}
  \begin{aligned}
   \frac{d}{dt}&\(S \d_{x}^{\beta } (x_k\bu^\e ),
  \d_{x}^{\beta } (x_k \bu^\e)\)  \le \\
&\le C \(\left\|\bu^\e\right\|_{H^s}+\left\|x\bu^\e\right\|_{H^{s-1}} 
  \)\(\|\bu^\e\|_{H^s}^2+ \|x \bu^\e\|_{H^{s-1}}^2\).  
  \end{aligned}
\end{equation}
The term $\|\bu^\e\|_{H^s}^2$ is due to the commutator $[x_k,L]$:
\begin{align*}
 \left|\( S \d^\beta_x \( [x_k,L]\bu^\e\),\d^\beta_x\(x_k
   \bu^\e\)\)\right| & = 
\left|\(  \d^\beta_x \( [x_k,L]\bu^\e\),\d^\beta_x\(x_k
   \bu^\e\)\)\right|\\ 
&\le \|\bu^\e\|_{H^s}\|x \bu^\e\|_{H^{s-1}}.
\end{align*}
Summing over the inequalities \eqref{eq:10h45} and \eqref{eq:10h46}
yields a closed set of estimates, from which we infer the analogue of
Proposition~\ref{prop:p};  note that the time $T_s$ is not larger
than $T$ by construction, and may be strictly smaller than $T$, due to
possible shocks for \eqref{eq:systhypV}. We also
mention the fact that 
the above analysis gives $v^\e = \nabla \phi^\e \in C([0,T_s];H^s)$,
with $x \nabla \phi^\e \in C([0,T_s];H^{s-1})$: back to \eqref{eq:modif0},
this shows that $ \phi^\e \in C([0,T_*];H^{s-1})$, but that we cannot
claim that $ x\phi^\e \in C([0,T_*];H^{s-1})$. 
\smallbreak

Passing to the limit $\e \to 0$ in \eqref{eq:modif0}, it is natural to
introduce the system:
\begin{equation}\label{eq:11h09}
\begin{aligned}
    \d_t \phi +\frac{1}{2}\left|\nabla \phi\right|^2 +\nabla
    \eik\cdot \nabla \phi+ f\(
    |a|^2\)&= 0,\\
    \d_t a +\nabla \phi \cdot \nabla a +\nabla \eik \cdot \nabla
    a +\frac{1}{2}a 
\Delta \phi +\frac{1}{2}a 
\Delta \eik &= 0,\\
\phi\big|_{t=0}=0\quad ;\quad a\big|_{t=0}&= a_0\, .
\end{aligned} 
\end{equation}
The above analysis shows that this system has a unique solution in
$H^s$ with the first momentum in $H^{s-1}$, 
locally in time for $t\in [0,T_*]$, for some $T_*\in ]0,T]$;  $T_*$ is
independent of $s$, from the usual continuation
principle, explained for instance in \cite[Section~2.2]{Majda}. We
easily obtain the analogue of 
Proposition~\ref{prop:estprec}: mimicking the above computations, we
can estimate the error $(a^\e -a,\phi^\e-\phi)$ in $H^s$, by first
estimating $(a^\e -a,\nabla\phi^\e-\nabla\phi)$ \emph{and its first
  momentum} in $H^s$. We deduce that $T_s\ge T_*$, and
\eqref{eq:BKWVtpetit} follows. 
\smallbreak

The end of Theorem~\ref{theo:BKWV} can then be proved as in
\cite{Grenier98}: from the above analysis, the functions
\begin{equation*}
  \frac{a^\e -a}{\e}\virgp\ \frac{\nabla \phi^\e -\nabla\phi}{\e}
\end{equation*}
and their first momentum, are bounded in $H^s$ for every $s\ge 0$. A
subsequence converges to the linearization of 
\eqref{eq:systhypV}, yielding a pair $(a^{(1)},\phi^{(1)})$. By
uniqueness for the limit system, the whole sequence is convergent, and
the analogue of Proposition~\ref{prop:correc} follows. This completes
the proof of Theorem~\ref{theo:BKWV}.


\section{Extension to the case $0<\kappa <1$}
\label{sec:kappaint}

When $0<\kappa <1$, we propose an analysis which can be considered as
a generalization of the study led in
Section~\ref{sec:surcrit}. Throughout this paragraph, we suppose that
Assumptions~\ref{hyp:geom}-\ref{hyp:surcrit} are satisfied. Again,
we write the exact solution as
\begin{equation*}
  u^\e = a^\e e^{i\Phi^\e/\e},\quad \text{with }\Phi^\e=\eik +\phi^\e.
\end{equation*}
The unknown function is the pair $(a^\e,\phi^\e)$. We have two unknown
functions to solve a single equation, \eqref{eq:nlssemi}. We can
choose how to balance the terms: we resume the approach followed when
$\kappa=0$. Note that this approach would also be efficient for the
case $\kappa \ge 1$, with the serious drawback that we still assume
$f'>0$, an assumption proven to be unnecessary when $\kappa\ge 1$ (see
Section~\ref{sec:sub}). We impose: 
\begin{equation}\label{eq:modifint}
\begin{aligned}
    \d_t \phi^\e +\frac{1}{2}\left|\nabla \phi^\e\right|^2 +\nabla
    \eik\cdot \nabla \phi^\e+ \e^{\kappa}f\(
    |a^\e|^2\)&= 0,\\
    \d_t a^\e +\nabla \phi^\e \cdot \nabla a^\e +\nabla \eik \cdot \nabla
    a^\e +\frac{1}{2}a^\e 
\Delta \phi^\e +\frac{1}{2}a^\e 
\Delta \eik &= i\frac{\e}{2}\Delta a^\e,\\
\phi^\e\big|_{t=0}=0\quad ;\quad a^\e\big|_{t=0}&= a^\e_0\, .
\end{aligned} 
\end{equation}
This is the same system as \eqref{eq:modif0}, with only $f$ replaced
by $\e^{\kappa}f$. Mimicking the analysis of
Section~\ref{sec:surcrit}, we work with the unknown $\bu^\e$ given by
the same definition: it solves the system \eqref{eq:systhypV}, where
only the matrices $A_j$  have changed, and now depend on $\e$. The
symmetrizer is the same as before, with $f'$ replaced by $\e^\kappa
f'$: the matrix $S=S^\e$ is not bounded as $\e \to 0$, but its inverse
is. We see that \eqref{eq:estdtS} and \eqref{eq:dtS}
still hold, independent of $\kappa$. We claim that inequalities
similar to \eqref{eq:10h45} and \eqref{eq:10h46} hold:
\begin{align*}
  \frac{d}{dt}\(S^\e \d_{x}^{\alpha } \bu^\e ,
  \d_{x}^{\alpha } \bu^\e\)  \le C
  \(\left\|\bu^\e\right\|_{H^s}+\left\|x\bu^\e\right\|_{H^{s-1}} 
  \)&\sum_{|\gamma|\le s}\(S^\e \d_{x}^{\gamma } \bu^\e ,
  \d_{x}^{\gamma } \bu^\e\),\\
\frac{d}{dt}\(S^\e \d_{x}^{\beta } (x_k\bu^\e ),
  \d_{x}^{\beta } (x_k \bu^\e)\)  
\le C \(\left\|\bu^\e\right\|_{H^s}+\left\|x\bu^\e\right\|_{H^{s-1}} 
  \)&\Big( \sum_{|\gamma|\le s}\(S^\e \d_{x}^{\gamma } \bu^\e ,
  \d_{x}^{\gamma } \bu^\e\)\\
 + \sum_{{|\gamma|\le s-1}\atop{1\le j\le
  n}}(S^\e \d_{x}^{\gamma 
  }& (x_j\bu^\e) , 
  \d_{x}^{\gamma } (x_j\bu^\e))\Big),
\end{align*}
where the map $C(\cdot)$ is locally bounded, \emph{independent of
}$\e\in]0,1]$. The fact that such estimates remain valid, with this
dependence upon $\e$, stems essentially from the following reasons:
\begin{itemize}
\item The matrices $S^\e$ and $B_j$ are diagonal.
\item The matrix $M$ is block diagonal: the blocks correspond to the
  presence/absence of $\e$ in $S^\e$. 
\item The matrices $S^\e A_j^\e$ are independent of $\e\in]0,1]$.
\item The inverse of $S^\e$ is uniformly bounded on compacts, as $\e
  \to 0$. 
\end{itemize}
A continuity argument and Gronwall lemma then imply the analogue of
Proposition~\ref{prop:p}: $\bu^\e$ exists locally in time, with
$H^s$-norm uniformly bounded as $\e\to 0$. Note that since
$\phi^\e\big|_{t=0}=0$, we have: 
\begin{equation*}
  \(S^\e \d_{x}^{\alpha } \bu^\e ,
  \d_{x}^{\alpha } \bu^\e\)\big|_{t=0} = \O(1),
\end{equation*}
and we infer more precisely:
\begin{equation*}
  \| a^\e \|_{L^\infty ([0,T_*];H^s)} =\O(1) \quad ; \quad
    \left\| \nabla \phi^\e \right\|_{L^\infty ([0,T_*];H^s)} =\O\(\e^\kappa\). 
\end{equation*}
It seems natural to change unknown functions, and work with
$\widetilde \phi^\e = \e^{-\kappa}\phi^\e$ instead of $\phi^\e$. With
this, we somehow correct the shift in the cascade of equations caused by
the factor $\e^\kappa$ in front of the nonlinearity. Then
\eqref{eq:modifint} becomes: 
\begin{equation}\label{eq:modifinttilde}
\begin{aligned}
    \d_t \widetilde \phi^\e +\frac{\e^\kappa}{2}\left|\nabla \widetilde
      \phi^\e\right|^2 +\nabla 
    \eik\cdot \nabla \widetilde \phi^\e+ f\(
    |a^\e|^2\)&= 0,\\
    \d_t a^\e +\e^\kappa\nabla \widetilde \phi^\e \cdot \nabla a^\e +\nabla
    \eik \cdot \nabla 
    a^\e +\frac{\e^\kappa}{2}a^\e 
\Delta \widetilde \phi^\e +\frac{1}{2}a^\e 
\Delta \eik &= i\frac{\e}{2}\Delta a^\e,\\
\widetilde \phi^\e\big|_{t=0}=0\quad ;\quad a^\e\big|_{t=0}&= a^\e_0\, .
\end{aligned} 
\end{equation}
The pair $(\widetilde \phi^\e, a^\e$) is bounded in
$C([0,T];H^s)$. Therefore,  a subsequence is convergent, and the limit
is given by: 
\begin{equation*}
\begin{aligned}
    \d_t \widetilde \phi  +\nabla 
    \eik\cdot \nabla \widetilde \phi+ f\(
    |a|^2\)&= 0\ ; \quad \widetilde \phi\big|_{t=0}=0, \\
    \d_t a +\nabla \eik \cdot \nabla a +\frac{1}{2}a 
\Delta \eik &= 0\ ; \quad a\big|_{t=0}= a_0\, .
\end{aligned} 
\end{equation*}
We see that $a$ solves \eqref{eq:alibre};  $\widetilde \phi $ is
given by an 
ordinary differential equation along the rays associated to $\eik$,
with a source term showing nonlinear effect: $f\( |a|^2\)$. By
uniqueness, the whole sequence is convergent. Roughly speaking, we see that if 
\begin{equation*}
  \bw^\e = \,^t \(\nabla \( \widetilde \phi^\e-\widetilde \phi\),a^\e -a\), 
\end{equation*}
then Gronwall lemma yields:
\begin{equation*}
  \(S^\e \d^\alpha_x \bw^\e , \d^\alpha_x \bw^\e\) \le C \( \e +
  \e^\kappa\)\le 2 C\e^\kappa.
\end{equation*}
We infer:
\begin{prop}\label{prop:estprecint}
Let $s>2+n/2$. Then \eqref{eq:modifint} has a unique solution
$(a^\e,\phi^\e)\in C([0,T];H^s)^2$, such that $x_k a^\e, x_k\d_j \phi^\e\in
C([0,T];H^s)$,  for every $1\le j,k\le n$ ($T$ is given by
Lemma~\ref{lem:hj}). Moreover, there exists 
$C_s$ independent of $\e$ such that for every $0\le t\le T$, 
\begin{equation}\label{eq:tpetitint}
  \| a^\e (t)- a(t)\|_{H^s}\le C_s \e^\kappa \quad ;\quad 
  \| \phi^\e(t) -\e^\kappa \widetilde \phi\|_{H^s} \le C_s \e^{2\kappa} t , 
\end{equation}
where $a$ is given by \eqref{eq:alibre}. 
\end{prop}

Three cases must be distinguished:
\begin{itemize}
\item If $1/2< \kappa<1$, then we can infer the analogue of
  \eqref{eq:BKWVtgrand}. 
\item If $\kappa =1/2$, then we can infer the analogue of
  \eqref{eq:BKWVtpetit} (but not yet of \eqref{eq:BKWVtgrand}).
\item If $0<\kappa <1/2$, then we must pursue the analysis, and
  compute a corrector of order $\e^{2\kappa}$.  
\end{itemize}
We shall not go further into detailed computations, but instead, discuss
the whole analysis in a rather loose fashion. However, we note that
all the ingredients have been given for a complete justification. 
\smallbreak

Let $N= [1/\kappa]$, where $[r]$ is the largest integer not larger than
$r>0$. We construct $a^{(1)},\ldots,a^{(N)}$ and
$\widetilde\phi^{(1)},\ldots,\widetilde\phi^{(N)}$ such that:
\begin{align*}
 & \left\| a^\e -a-\e^\kappa a^{(1)}-\ldots -\e^{N\kappa}a^{(N)}
  \right\|_{L^\infty([0,T];H^s)} +\\
&+\left\| \widetilde\phi^\e -\widetilde \phi -\e^{\kappa} 
  \widetilde \phi^{(1)}- \ldots -\e^{N\kappa}\widetilde \phi^{(N)}
  \right\|_{L^\infty([0,T];H^s)}  
  =o\( \e^{N\kappa}\).  
\end{align*}
But since $N +1 >1/\kappa$, we have:
\begin{equation*}
  \left\| \phi^\e -\e^{\kappa} \widetilde \phi -\e^{2\kappa} 
  \widetilde \phi^{(1)}- \ldots -\e^{N\kappa}\widetilde \phi^{(N-1)}
  \right\|_{L^\infty([0,T];H^s)} 
  =\O\( \e^{(N+1)\kappa}\)=o(\e).
\end{equation*}
The analogue of \eqref{eq:BKWVtgrand} follows:
\begin{equation*}
  \left\| u^\e - a e^{ i \eik/\e + i\phi_{\rm app}^\e}
  \right\|_{L^\infty([0,T];L^2\cap  L^\infty)} =o(1),
\end{equation*}
where 
\begin{equation*}
  \phi_{\rm app}^\e = \frac{\widetilde \phi}{\e^{1-\kappa}} +\frac{\widetilde
  \phi^{(1)}}{\e^{1-2\kappa}}+\ldots + 
  \frac{\widetilde \phi^{(N-1)}}{\e^{1-N\kappa}}.
\end{equation*}
\begin{rema}
  In the case $\kappa =1$, $N=1$, and the above analysis shows that one phase
  shift factor appears: we retrieve the function $G$ of
  Proposition~\ref{prop:sub} (under the unnecessary assumption $f'>0$). If
  $\kappa >1$, then $N=0$, and we see that $a e^{i\eik/\e}$ is a good
  approximation for $u^\e$. 
\end{rema}



\bibliographystyle{amsplain}
\bibliography{../../carles}

\end{document}